\documentclass[11pt,twoside]{amsart}
\usepackage{amsmath}
\usepackage{amsfonts}
\usepackage{amssymb}
\usepackage{amsthm}
\usepackage[english]{babel}

\newtheorem{theo}{Theorem}
\newtheorem{lem}[theo]{Lemma}
\newtheorem{prop}[theo]{Proposition}

\newcommand{\E}{\ensuremath{\mathbb {E}}}
\renewcommand{\P}{\ensuremath{\mathbb {P}}}

\newcommand{\N}{{\mathbb N}}
\newcommand{\G}{{\mathcal G}}

\newcommand{\F}{{\mathcal F}}

\newcommand{\B}{{\cal B}}

\newcommand{\M}{{\mathcal M}}

\def\E{\mathbb{E}}
\def\Z{\mathbb{Z}}
\def\B{{\mathcal B}}
\def\F{{\mathcal F}}
\def\R{\mathbb{R}}

\begin{document}

\title
[Limit theorems for $GL_d(\R)$]
{ Limit theorems for the left random walk on $GL_d(\R)$}


\author{Christophe Cuny }
\address{Laboratoire MICS, Centrale-Supelec, Grande Voie des Vignes, 92295 Chatenay-Malabry cedex, France.}
\email{christophe.cuny@centralesupelec.fr}

\author{J\'er\^ome Dedecker}
\address{Laboratoire MAP5 (UMR 8145), Universit\'e Paris Descartes, Sorbonne Paris Cit\'e,  
45 rue des Saints P\`eres, 75270 Paris Cedex 06, France}
\email{jerome.dedecker@parisdescartes.fr}

\author{Christophe Jan}
\address{Lyc\'ee Claude Fauriel, Avenue de la lib\'eration, 42000 Saint-Etienne, France}
\email{christophe.jan@ac-lyon.fr}
\thanks{}


\begin{abstract}
Motivated by a recent work of Benoist and Quint and extending results 
from the PhD thesis of the third author, we obtain limit theorems 
for products of independent and identically distributed elements  of $GL_d(\R)$, such as 
the Marcinkiewicz-Zygmund strong law of large numbers, the CLT (with rates 
in Wasserstein's distances) and almost sure invariance principles with rates. 
\end{abstract}

\maketitle

\section{Introduction}

Let $(Y_n)_{n\ge 1}$ be independent random matrices taking values 
in $G:=GL_d(\R)$, $d\ge 2$ (the group of invertible $d$-dimensional real matrices),
with common distribution $\mu$. Let $\|\cdot\|$ be the euclidean norm on $\R^d$. We wish to study the asymptotic behaviour 
of $(\log \|Y_n \cdots Y_1\|)_{n\ge 1}$, where for every 
$g\in GL_d(\R)$, $\|g\|:=\sup_{x, \|x\|=1}\|gx\|$. 

\smallskip 

We shall say that $\mu$ has a (polynomial)  moment of order $p\ge 1$, if 
\begin{equation}\label{moment-p}
\int_G (\log N(g))^p \mu(dg)<\infty \, ,
\end{equation}
where $N(g):=\max(\|g\|,\|g^{-1}\|)$.

It follows from Furstenberg and Kesten \cite{FK} that, as soon as 
$\mu$ admits a moment of order 1,
$$
\lim_{n\to +\infty} \frac1n \log \|Y_n \cdots Y_1\|= \lambda_\mu 
\qquad \mbox{$\P$-a.s.}\, , 
$$
where $\lambda_\mu:= \lim_{n\to +\infty} n^{-1}
\E(\log \|Y_n \cdots Y_1\|)$ is the so-called first Lyapounov exponent. 

\smallskip

If moreover, no proper subspace of $\R^d$ is invariant by the closed 
semi-group generated by the support of $\mu$, then (see for instance 
Proposition 7.2 page 72 in \cite{BL}), for every $x\in \R^d-\{0\}$, 
\begin{equation}\label{SLLN}
\lim_{n\to +\infty} \frac1n \log \|Y_n \cdots Y_1x\|= \lambda_\mu 
\qquad \mbox{$\P$-a.s.}\, , 
\end{equation}

\smallskip

Our goal is to study the rate in the above convergences, assuming 
higher moments (and stronger algebraic conditions), as well as the Central Limit Theorem (CLT) or the 
Law of the Iterated Logarithm (LIL), and the rates of convergence in those limit theorems. 

\smallskip

The CLT question benefited from several papers 
under an exponential moment, i.e. 
$\int_G (N(g))^\alpha \mu(dg) <\infty$ for some $\alpha>0$, and 
some algebraic conditions, see the next section for more details. Let us 
mention among others the papers by Le Page \cite{Lepage} and
Guivarc'h and Raugi \cite{GR}.

\smallskip

Quite recently, Benoist and Quint \cite{BQ} proved the CLT 
under the existence of a moment of order 2. Their proof is based 
on Gordin's martingale approximation method. By an elegant but somewhat tricky argument, they provide an explicit martingale-coboundary decomposition adapted to the problem. Moreover, as intermediary steps,
they proved a result about \emph{complete convergence} as well as an 
integrability property with respect to the invariant 
probability measure on $X:=P_{d-1}(\R)$ (the projective space of $\R^d$), see the next section for further details and definitions. 
Let us  mention here  that most of the results of \cite{BQ} hold for linear groups on any local field. 

\smallskip

Rates in the CLT under polynomial moments have been announced in 
Jan \cite{Jan-CRAS} (with proof in \cite{Jan}) and the CLT has been 
proved in the PhD thesis of the third author \cite{Jan} under a 
moment of order $2+\varepsilon$, for any $\varepsilon >0$. 
His method of proof is also based on martingale approximation, 
but relies on estimates that seem more suitable to obtain precise rates of convergence (in the CLT and the strong invariance principle) than 
the approach of Benoist and Quint, at least in the case of $GL_d(\R)$.

\smallskip In Section \ref{Sec:Results} below, we give our main 
results for the sequence $(\log \|Y_n \cdots Y_1 x \|)_{n \geq 0}$
and any starting point $x \in {\mathbb R}^d- \{0\}$. 
We follow the approach described in Jan's PhD thesis \cite{Jan} (refining some of his computations), combined with recent or new results 
about rates in the strong invariance principle  and rates in the CLT 
(see Section 3).
At the very end of the paper (cf. Section  \ref{Sec:MNMC}), we also borrow one main argument from Benoist and Quint \cite{BQ}, to prove that the rates
of convergence in the CLT apply  to the sequence $(\log\|Y_n \cdots Y_1\|)_{n\ge 1}$,  and to obtain some results 
for  the sequence of matrix coefficients $( \log |\langle Y_n \cdots Y_1 x,y\rangle|)_{n\ge 1}$. 
In the same final section, we also briefly explain how to weaken the assumption of proximality (see the next section for the definition) by 
using another argument from \cite{BQ}.


\section{Results} \label{Sec:Results}

Let $G:=GL_d(\R)$, $d\ge 2$, endowed with its Borel $\sigma$-algebra 
$\B(G)$. Let $X:=P_{d-1}(\R)$ be the projective space 
of ${\mathbb R}^d-\{0\}$, and write $\overline x$ as the projection of $x \in {\mathbb R}^d- \{0\}$ to $X$. 
Then $G$ acts continuously on  $X$ in a natural way : $g\cdot \overline x=  \overline{gx}$.

\smallskip

Let $\mu$ be a probability measure  on $\B(G)$. Denote by $\Gamma_\mu$ 
the closed semi-group generated by the support of $\mu$. Assume that 
$\mu$ is \emph{strongly irreducible}, i.e. that no proper finite union of subspaces of $\R^d$ are invariant by $\Gamma_\mu$ and that 
it is \emph{proximal}, i.e. that there exists a matrix in $\Gamma_\mu$ 
admiting a unique (with multiplicity one) eigenvalue with maximum modulus. 

\smallskip

For such a measure $\mu$, 
it is known  that there exists a unique invariant measure 
$\nu$ on $\B(X)$ (see for instance Theorem 3.1 of \cite{BL}) in the following sense:  for any continuous
and bounded function $h$ from $X$ to ${\mathbb R}$
$$
  \int_X h(x) \nu(dx) = \int_G \int_X  h(g \cdot  x) \mu(dg) \nu(dx) \, .
$$

\smallskip

We consider the left random 
walk of law $\mu$ on $X$. 
Let us recall its construction. 

\smallskip 

Let $\Omega:=X\times G^{\N^*}$ and $\F:=\B(X)\otimes \B(G)^{\otimes \N^*}$, 
where $\N^*=\{1,2,\ldots\}$. 
For every  probability measure $\tau$ on $\B(X)$,  we define $\P_\tau:=\tau \otimes \mu^{\otimes \N^*}$. As usual we note $\P_{\overline x}
:=\P_{\delta_{\overline x}}$, for every $\overline x \in X$.  Define the coordinate process 
$(Y_n)_{n\in \N}$ ($\N=\{0,1,\ldots\}$), i.e. 
$Y_0((\overline x,g_1,g_2,\ldots))=\overline x$ and for every $n\in \N^*$, 
$Y_n((\overline x,g_1,g_2,\ldots))=g_n$, and then $\F_n$, the $\sigma$-algebra generated by $\{Y_0,\ldots ,Y_n\}$.

\smallskip

 Finally, define a measurable transformation 
$\eta$ on $\Omega$ by 
$$
\eta((\overline x,g_1,g_2,\ldots))=
(g_1\cdot \overline x ,g_2,g_3,\ldots) \, .
$$ 

\smallskip

The left random walk of law $\mu$ is the process $(W_n)_{n\in \N}$, defined by $W_0:= Y_0$ and for every $n\in \N^*$, $W_n=W_0\circ \eta^n$. Hence,  it is a Markov
chain defined by the recursive equation
$W_n=Y_n W_{n-1}$ for $n \in \N^*$. 

\smallskip

Recall that for every probability measure 
$\P_\tau$, $Y_0$ is a random variable with law $\tau$ independent from 
the  sequence $(Y_n)_{n\in \N^*}$ of independent and identically distributed (iid) random variables. 
Recall also that $\P_\nu$ is $\eta$-invariant hence, under $\P_\nu$, $(W_n)_{n\in \N}$ is identically distributed with 
common marginal distribution $\nu$. Moreover, since $\nu$ is the unique $\mu$-invariant probability, then $(\Omega,\F,\P_\nu,\eta)$ is ergodic 
(see e.g. Proposition 1.14 page 36 of \cite{BQ-book}).

\smallskip

We want to study the process $(X_n)_{n\in \N^*}$ given by 
$X_n:=\sigma(Y_n,W_{n-1})$ for every $n\in \N^*$, where 
for every $g\in G$ and every $\bar x \in X$, 
$$
\sigma(g,\bar x)= 
\log \left (\frac{\|g\cdot x\|}{\|x\|} \right )\, .
$$

\smallskip

Let us denote $A_0={\mathrm {Id}}$ and, for every $n\in \N^*$, $A_n:= Y_n\cdots Y_1$, so that $X_n=\sigma(Y_n,A_{n-1}W_0)$.
Let $S_n:= X_1+\cdots +X_n$, and note that 
$S_n= \log \|A_n W_0^*\|$, where $W_0^*$ is 
an element of ${\mathbb R}^d$ such that 
$\overline{ W_0^*}= W_0$ and $\|W_0^*\|=1$.
Finally, let
$$
B_n=
\left \{ \frac{S_{[nt]}-[nt] \lambda_{\mu}}{\sqrt n}-\frac{(nt-{[nt]})}{\sqrt n} (X_{[nt]+1}-  \lambda_{\mu}), t \in [0,1] \right \}
$$
be the partial sum process with values in the space 
$C([0,1])$ of continuous functions on $[0,1]$ equipped with the uniform metric.

As usual in the Markov chain setting, we  denote 
by $X_{n, \overline x}$ the random variable 
$X_n$ 
for which $W_0= \overline x$.
Let also
$S_{n, \overline x}$ be the corresponding partial sum, 
and  $B_{n, \overline x}$ be the corresponding 
process. Note that 
$S_{n, \overline x}= \log \|A_n x \|$ if $\|x\|=1$.

Note that the distribution of  the sequence $(X_{n, \overline x})_{n\in \N^*}$ is the same for any 
probability  $\P_\tau$ on $\Omega$ (in fact
$(X_{n, \overline x})_{n\in \N^*}$ is a function 
of $\overline x$ and $(Y_n)_{n\in \N^*}$, so that its distribution depends 
only on $\overline x$ and $\mu$). 
Hence, we shall write
``${\mathbb P}$-almost surely" (${\mathbb P}$-a.s.)  instead of ``${\mathbb P}_{\tau}$-almost surely'', and
${\mathbb E}(\cdot)$ instead of 
${\mathbb E}_{\tau}(\cdot)$, for all the quantities 
involving the sequence $(X_{n, \overline x})_{n\in \N^*}$
(and more generally for all the quantities involving only 
the sequence $(Y_n)_{n\in \N^*}$).
With these notations,  for any positive and  measurable function
$f$, 
$$\E(f(X_{n,\overline x}))= \E_{\overline x}(f(X_n))=
\E (f(X_n)|W_0=\overline x) \, .
$$

\smallskip

Our study will only require polynomial moments for $\mu$. 
As already mentionned, when $\mu$ has a moment of order 1, the strong law of large numbers \eqref{SLLN} holds for any starting point. Moreover, one can identify the limit $\lambda_\mu$ $via$ the ergodic theorem for strictly stationary sequences. 
It follows that,  for every $\overline x\in X$, 
$$\frac{S_{n, \overline x}}{n}\underset{n\to +\infty}\longrightarrow \lambda_\mu= \int_G\int_X \sigma(g,u) \mu(dg)\nu(du)\quad \text{${\mathbb P}$-a.s.},$$ 
 see for instance Corollary 3.4 page 54 of \cite{BL} or 
Theorem 3.28 of \cite{BQ-book}. Our goal is to strengthen that strong law of large numbers 
when higher moments are assumed.

As already mentionned, in the next theorem, item $(ii)$ has been obtained by Benoist and Quint \cite{BQ}. As observed in the introduction of \cite{BQ}, 
their method also allow to prove item $(ii)$ of the next theorem 
when $p=2$.

\begin{theo}\label{main-theo}
Let $\mu$ be a proximal and strongly irreducible probability measure on 
$\B(G)$. Assume that $\mu$ has a moment of order $p\ge 1$. 
\begin{itemize}
\item [$(i)$] If $1\le p<2$ then, for every $ \overline x\in X$, 
$$\frac{S_{n, \overline x}-n\lambda_\mu}{ n^{1/p}}
\underset{n\to +\infty}\longrightarrow 0 \quad  
\text{${\mathbb P}$-a.s.}.$$
\item [$(ii)$] If $p=2$ then $ n^{-1}\E_\nu ((S_n-n\lambda_\mu)^2)
\longrightarrow \sigma^2$ as $n \rightarrow \infty$, and, for any continuous and bounded function
$\varphi$ from $C([0,1])$ (equipped with the sup norm) to ${\mathbb R}$, 
$$ 
 \lim_{n \rightarrow \infty} \sup_{\overline x \in X} \left | 
 {\mathbb E}( \varphi(B_{n, \overline x})) - 
 \int \varphi ( \sigma \varpi) w(d\varpi)
 \right  | = 0 \, ,
$$
where $w$ is the distribution of a standard Wiener process.
\item [$(iii)$] If $2\le p< 4$ then, for every (fixed) $\bar x 
\in X$, one can redefine $(S_{n,\bar x})_{n\ge 1}$ 
without changing its distribution on a (richer) probability space  on which there exists iid random variables 
$(W_n)_{n\ge 1}$ with common distribution ${\mathcal N}(0,\sigma^2)$,  such that, 
$$
\left |S_{n,\bar x}-n\lambda_\mu- \sum_{i=1}^n
W_i \right |= o(r_n)
\quad \text{${\mathbb P}$-a.s.}\, ,
$$
where $r_n=\sqrt{n\log \log n}$ when  $p=2$ and  $r_n= n^{1/p} \sqrt{\log n}$ 
when $2<p<4$.
\item [$(iv)$] If $p=4$ then, for \emph{every} (fixed) $\bar x 
\in X$, one can redefine $(S_{n,\bar x})_{n\ge 1}$ 
without changing its distribution on a (richer) probability space  on which there exists iid random variables 
$(W_n)_{n\ge 1}$ with common distribution ${\mathcal N}(0,\sigma^2)$,  such that, 
$$\left |S_{n,\bar x}-n\lambda_\mu- \sum_{i=1}^n 
W_i
\right |= O\left (n^{1/4}\sqrt{\log n} \left(\log \log n \right)^{1/4}  \right)
\quad \text{${\mathbb P}$-a.s.}\, .
$$
\end{itemize}
\end{theo}
\noindent {\bf Remark. } 
Let us recall the famous result by  Koml\'os, Major and Tusn\'ady \cite{KMT}. Let $(V_n)_{n \in \N}$ be iid 
variables in ${\mathbb L}^p$, $p>2$. Then, extending the 
probability space if necessary, it is possible to construct iid  random variables 
$(Z_n)_{n\ge 1}$ 
with common distribution ${\mathcal N}(0,{\mathrm{Var}} (X_1))$ such that
$$\left | \sum_{i=1}^n 
(V_i-\E(V_i))- \sum_{i=1}^n 
Z_i
\right |= o\left (n^{1/p}  \right)
\quad \text{a.s.}\, .
$$
Hence, for $p \in (2,4]$, our results are close to the iid
situation. The logarithmic loss seems to be difficult to avoid with our
approach based on  martingale approximation.

\smallskip
 
 \noindent {\bf Remark. } It follows from Theorem 4.11 $c)$ of 
 \cite{BQ} that $\sigma\neq 0$ when $\Gamma_\mu$ has unbounded image 
 in  $PGL(V)$. 

\medskip

The proof of Theorem \ref{main-theo} will result from general limit theorems under projective 
conditions. When $1<p<2$ those results are new and when $p>2$, 
the obtained rates slightly improve previous results (see for instance \cite{DDM}). 

\medskip

We also obtain rates of convergence for Wasserstein's distances in the 
central limit theorem. Let us first recall the definition of these minimal distances. 
Let ${\mathcal L}(\nu_1, \nu_2)$ be the set of the probability laws on
$\mathbb R^2$ with  marginals $\nu_1$ and $\nu_2$. The Wasserstein distances of 
order $r$ between $\nu_1$ and $\nu_2$
are defined as follows: 
$$
W_r(\nu_1, \nu_2) =  \left\{
\begin{array}{ll} \displaystyle \inf \left \{ \int |x-y|^r P(dx, dy)  :  P
\in {\mathcal L}(\mu, \nu) \right \} & \text{if $0< r < 1$} \\
\displaystyle \inf \left \{ \left (\int |x-y|^r P(dx, dy) \right )^{1/r}
: P \in {\mathcal L}(\mu, \nu) \right \}  & \text{if $r \geq 1$} \, .
\end{array}
\right.
$$
It is well known that, for $r \in (0, 1]$,
$$
W_{r} ( \nu_1 , \nu_2)   =  \sup \left \{  \nu_1 (f) -  \nu_2 (f)
: f \in \Lambda_r \right \},
$$
where $\Lambda_r$ is the set of $r$-H\"older functions such that
$
|f(x)-f(y)|\leq |x-y|^r
$
for any reals $x,y$.
For $r \geq 1$, one has
\begin{equation*} \label{def2wasser}
W_r(\nu_1, \nu_2) =  \left ( \int_0^1 |F_1^{-1}(u) - F_2^{-1}(u) |^r du\right
)^{1/r},
\end{equation*}
where $F_1$ and $F_2$ are the respective distribution functions of $\nu_1$ and $\nu_2$, and $F_1^{-1}$ and $F_2^{-1}$ 
are their generalized inverse.

We obtain 

\begin{theo}\label{rates}
Let $\mu$ be a proximal and strongly irreducible probability measure on 
$\B(G)$. For any $\overline x  \in X$, denote by $ \nu_{n,\overline x}$ 
 the distribution of $n^{-1/2}(S_{n, \overline x} - n \lambda_\mu)$.  Let also
$G_{\sigma}$  be the normal distribution  with mean zero and variance $\sigma^2$
given in Theorem  \ref{main-theo}$(ii)$ (provided 
$\mu$ has a moment of ordrer 2).
\begin{itemize}
\item  [$(i)$] Assume that $\mu$ has a moment of order $p\in (2,3)$. Then, for any  $r \in [p-2, p]$, 
$$
\sup_{\overline x\in X} W_r \left( \nu_{n, \overline x},G_{\sigma} \right)= O\left (n^{-(p-2)/ 2 \max (1, r)} \right ) \, .
$$
\item [$(ii)$] Assume that $\mu$ has a moment of order 3. Then, for any $r$ in $ (1, 3]$, 
$$
\sup_{\overline x\in X} W_r \left(\nu_{n, \overline x}, G_{\sigma} \right)= O\left (n^{-1/ 2 r} \right ), 
$$
and for $r=1$, 
\begin{equation}\label{rate1,3}
\sup_{\overline x\in X} W_1 \left(\nu_{n,\overline x},G_{\sigma} \right)
= O\left (n^{-1/ 2} \log n \right ) \, .
\end{equation}
\end{itemize}
\end{theo}

\noindent {\bf Remark. }
Except for 
$p=3, r=1$, the rates given in Theorem \ref{rates} are consistent 
with the iid case, in the following sense: let 
$(V_i)_{i \geq 1}$ be a sequence of iid random variables, 
where the $V_i$'s are centered and have a moment of order $p \in (2,3)$. 
Let $\nu_n$ be the distribution of 
$n^{-1/2}(V_1 + \cdots + V_n)$.
Then the rates given in Theorem \ref{rates} hold 
for $\nu_n$ instead of $\nu_{n,\overline x}$ and 
$\sigma^2= {\mathbb E}(V_1^2)$. Moreover, these
are the best known rates under the stated conditions (see 
the introduction of the paper \cite{DMR}). For
$p=3, r=1$ the rate in the iid case is $O(n^{-1/2})$, 
so there is a loss of order $ \log n$ in \eqref{rate1,3}.

\smallskip

\noindent {\bf Remark. }
Starting from Remark 2.3 of \cite{DMR}, we derive from
Theorem \ref{rates}   the following rates of convergence in the Berry-Esseen theorem:
If $\mu$ has a moment of order $p\in (2,3)$, then
$$
\sup_{ \overline x \in X} \sup_{t\in \R} \left |
\P\left (n^{-1/2}(S_{n, \overline x}-n\lambda_\mu)\le t \right )-\phi_\sigma(t) \right |
\leq O \left (  n^{-(p-2)/2(p-1)}\right ) \, ,
$$
where $\phi_\sigma$ is the distribution function 
of $G_{\sigma}$.
If $\mu$ has a moment of order $3$, then
$$
\sup_{ \overline x \in X} \sup_{t\in \R} \left |
\P\left (n^{-1/2}(S_{n, \overline x}-n\lambda_\mu)\le t \right )-\phi_\sigma(t) \right |
\leq O \left (  n^{-1/4} \sqrt {\log n }\right ) \, .
$$
Note that, when $\mu$ has moments of any order, 
Jan \cite{Jan}  obtained the rate $O(n^{-a})$ for any $a < 1/2$
in the Berry-Esseen theorem. 



\section{Auxiliary results on the cocycle}

In all this section $\mu$ is a proximal and  strongly irreducible 
probability measure on $\B(G)$. 
Let
 $\tilde X_{k} = X_{k} - \lambda_\mu$ and 
 $\tilde X_{k, \overline x} = X_{k, \overline x} - \lambda_\mu$. For $p \geq 1$, let $\|\cdot \|_{p, \tau}$ be the 
 ${\mathbb L}^p$-norm with respect to the probability
 $\P_\tau$ on $\Omega$. For the quantities involving 
 $ X_{k, \overline x}$, we shall write $\|\cdot \|_{p}$
 instead of $\|\cdot \|_{p, \tau}$, in accordance with the notations of Section \ref{Sec:Results}.

\smallskip

The proofs of Theorem \ref{main-theo} and  Theorem \ref{rates} will make use 
of general results for stationary sequences under projective conditions, 
i.e. conditions relying on the quantities 
$\|\E(\tilde X_n|W_0)\|_{p, \nu}$ for $p\geq 1$ and $\|\E(\tilde X_n \tilde X_k|W_0)-\E_\nu(\tilde X_n \tilde X_k)\|_{p/2, \nu}$
for $p \geq 2$. 

\smallskip

Those quantities were already studied in \cite{Jan}, where polynomial rates of convergence (to 0) were obtained. By refining the arguments 
of \cite{Jan} we  obtain the following improvements.

\begin{prop}\label{newprop}
 Assume that $\mu$ has a moment of order $p>1$.
Then, for $q \in [1,p)$,  
\begin{equation}\label{decadix}
\sum_{k=1}^\infty k^{p-q-1} 
\sup_{\overline x,  \overline y \in X} 
\E \left (|X_{k,\overline x}-X_{k, \overline y}|^{q} \right)< \infty \, \, .
\end{equation}
and for $q \in (0,1]$, 
\begin{equation}\label{decadix2}
\sum_{k=1}^\infty k^{p-2} 
\sup_{\overline x,  \overline y \in X} 
\E \left (|X_{k,\overline x}-X_{k, \overline y}|^{q} \right)  < \infty \, \, .
\end{equation}
\end{prop}

\noindent{\bf Remark.}
Since $\E(X_{k, \overline x})=\E_{\overline x}(X_k)=
\E(X_k | W_0= \overline x)$, and since $\E_\nu(X_k)= \lambda_\mu$, we easily infer from \eqref{decadix} that
\begin{equation}\label{simple-series}
 \sum_{k\ge 1} k^{p-2}\sup_{\overline x \in X}|\E_{\overline x}(X_k)-\lambda_\mu| <\infty \, .
 \end{equation}
In particular, using that $p+1/p>2$ whenever 
 $p>1$, it follows from  \eqref{simple-series} that 
\begin{equation}\label{borneMZ}
 \sum_{k\ge 1} k^{-1/p}\sup_{\overline x \in X}|\E_{\overline x}(X_k)-\lambda_\mu| <\infty\, .
\end{equation}
\noindent {\bf Remark.} Let us notice that the third author 
\cite{Jan} proved that for every $p\ge 2$ 
and every $\alpha\in [0,1)$, there exists $C_{p,\alpha}$ such that 
$\sup_{\overline x,\overline y\in X} \E(|X_{k,\overline x}
-X_{k,\overline y}|)\le \frac{C_{p,\alpha}}{k^{p(\alpha-1/2)}}$. In particular, when $p=2$, using a Theorem of Maxwell and Woodroofe \cite{MW}, 
that estimate is sufficient for the CLT   under a second moment and 
even for the invariance principle  (see Peligrad and Utev \cite{PU}). Hence,   the full conclusion of item $(ii)$ in Theorem \ref{main-theo} follows from the latter estimate. 
 
 \medskip

We shall also need the following controls.

 \begin{prop}\label{newprop2}
 Assume that $\mu$ has a moment of order $p>2$. Then
 \begin{equation}
 \label{square-series}\sum_{k\ge 1} k^{p-3}
 \sup_{\overline x , \overline y \in X}
 \E \left (\left|\tilde X_{k,\overline x}^2-\tilde X_{k, \overline y}^2 \right | \right)<\infty\, ,
 \end{equation}
 and for every $\gamma<p-3+1/p$,
 \begin{equation}
 \label{mixed-series}\sum_{k\ge 1} k^{\gamma}\sup_{\overline x , \overline y \in X}\sup_{ k
 \le j<i \le 2k}\E \left (
 \left |\tilde X_{i,\overline x}\tilde  X_{j,\overline x}-\tilde X_{i, \overline y}\tilde X_{j,\overline x} \right | \right)<\infty\, .
 \end{equation}
 \end{prop}
 \noindent {\bf Remark.}
 As in the previous remark, we easily infer that  
 \begin{equation}
 \label{square-series2}\sum_{k\ge 1} k^{p-3}
 \sup_{\overline x  \in X}
 \left | \E_{\overline x}\left (\tilde X_k^2 \right)-\E_\nu \left (\tilde X_k^2 \right )  \right | <\infty\, ,
 \end{equation}
 and for every $\gamma<p-3+1/p$,
 \begin{equation}
 \label{mixed-series2}\sum_{k\ge 1} k^{\gamma}
 \sup_{\overline x \in X}\sup_{ k
 \le j<i \le 2k}\left |\E_{\overline x}\left (\tilde X_i \tilde X_j \right )- \E_\nu\left (\tilde X_i \tilde X_j \right ) \right | <\infty\, .
 \end{equation}
 
 \medskip
 
 The proof of Propositions \ref{newprop} 
 and \ref{newprop2}  are  based on two auxiliary 
 lemmas. The first one gives  the regularity of 
the cocycle $\sigma$ with respect to a suitable metric, that we introduce right now.

\medskip

For $\overline x, \overline y \in X$, define 
$$
d(\overline x , \overline y):= \frac{\|x\wedge y\|}{\|x\|\, \|y\|}\, ,
$$ 
where $\wedge$ stands for the exterior product, see e.g. 
\cite[page  61]{BL} for the definition and some properties. Then, $d$ 
is a metric on $X$.

\medskip

For every $q>0$, define a non decreasing, concave function 
$H_q$ on $[0,1]$ by $H_q(0)=0$ and for every $x\in (0,1]$,
$H_q(x)= \big|\log (x{\rm e}^{-q-1})\big|^{-q}$. 

\smallskip


The next lemma may be seen as a version of Lemma 17 of Jan \cite{Jan}. 
\begin{lem}\label{lem-reg}
For every $\kappa>1$, there exists 
$C_\kappa>0$ such that for every $g\in G$ and every $\overline x, \overline y \in X$, 
\begin{equation}\label{regularity}
|\sigma(g,\overline x)-\sigma(g, \overline y)|\le C_\kappa (1+
\log N(g))^\kappa \,  H_{\kappa-1}(d(\overline x, \overline y))\, .
\end{equation}
\end{lem}
\noindent   {\bf Proof.} By Lemma 12.2 of \cite{BQ-book}, there exists $C>0$ such that for every $\overline x, \overline y\in X$,
 \begin{equation}\label{BQ-est}
 |\sigma(g,\overline x)- \sigma(g,\overline y)|\le C N(g)d(\overline x,\overline y)\, .
 \end{equation}

Now, it is not hard to prove that (notice that $\|g^{-1}\|^{-1}
\le  \|x\|^{-1}\|g x\|\le \|g\|$ for every $g\in G$ and every 
$x\in \R^d-\{0\}$), for every $\overline x\in X$ and every 
$ g\in G$,
\begin{equation}\label{trivial-est}
\sigma(g,\overline x)\le \log(N(g)) \, .
\end{equation}

 Assume that $d(\overline x, \overline y)\le 1/N(g)$. Using that 
$t \mapsto t(H_{\kappa-1}(t))^{-1}$ 
is non decreasing on  $(0, {\rm e}]$ and that   $N(g)\ge 1$, we have 
$$N(g) d(\overline x, \overline y)\le   
\frac{H_{\kappa-1}(d(\overline x,\overline y))}{H_{\kappa-1}(N(g)^{-1})}\, .
$$
Hence, by \eqref{BQ-est}, 
\begin{multline}\label{first-est}
|\sigma(g,\overline x)- \sigma(g,\overline y)|\le C H_{\kappa-1}(N(g)^{-1}) \, H_{\kappa-1}(d(\overline x,\overline y)) \\
\le C(\kappa +\log N(g))^{\kappa-1}H_{\kappa-1}(d(\overline x,\overline y))\, .
\end{multline}

Assume now that $ d(\overline x,\overline y)>1/N(g)$. By \eqref{trivial-est}, 

\begin{equation}\label{second-est}
|\sigma(g,\overline x)- \sigma(g,\overline y)|\le \frac{\log N(g)
 H_{\kappa-1}(N(g)^{-1})}
{H_{\kappa-1}(N(g)^{-1})}\le (\kappa+\log (N(g))^{\kappa}
H_{\kappa-1}(d(\overline x,\overline y))
\,.
\end{equation}

Combining \eqref{first-est} and \eqref{second-est}, we see that 
\eqref{regularity} holds. \hfill $\square$

\medskip

The next lemma is a result about complete convergence that may be derived 
from Proposition 4.1 of Benoist and Quint \cite{BQ}. A different proof is given in Section \ref{Sec:inter}.

\begin{lem}\label{lem-comp}
Assume that $\mu$ has a moment of order $p>1$. Then, there exists $\ell>0$, such that 
\begin{equation}\label{complete-cocycle}
\sum_{k\ge 1} k^{p-2}\max_{k\le j\le 2k}\sup_{\overline x,\overline y\in X, \overline x\neq \overline y} \P \left(
 \log \left(d(A_{j-1}\cdot \overline x, 
A_{j-1}\cdot \overline y)\right)\ge -\ell k\right)<\infty\, .
\end{equation}
\end{lem}

 \noindent {\bf Proof of Proposition \ref{newprop}.}  
 Let $\overline x , \overline y \in X$. Let $\ell>0$ be 
 as in Lemma \ref{lem-comp}. 
We start from the elementary inequality: If 
$A=\{\log N(Y_k)
\ge k\}$ and $B=\{\log d(A_{k-1}\cdot x, A_{k-1}\cdot y)\ge -\ell k\}$, 
\begin{multline}
\label{first}
|X_{k, \overline x}-X_{k, \overline y} | \leq 
|\sigma(Y_k, A_{k-1} \overline x) -\sigma(Y_k, A_{k-1} \overline y) |{\bf 1}_A 
+|\sigma(Y_k, A_{k-1} \overline x) -\sigma(Y_k, A_{k-1} \overline y) |{\bf 1}_B \\
+|\sigma(Y_k, A_{k-1} \overline x) -\sigma(Y_k, A_{k-1} \overline y) |{\bf 1}_{\{A^c \cap B^c \}}
\, .
\end{multline}

Using \eqref{trivial-est} and \eqref{regularity} (with $\kappa=(p+q)/q$), we infer 
from \eqref{first} that
\begin{multline*}
\E \left ( |X_{k, \overline x}-X_{k, \overline y} |^q \right )
   \le  C \E \left ( |\log N(Y_k)|^q {\bf 1}_{A} \right ) 
+ C \E \left (\left |\log N(Y_k)\right |^q {\bf 1}_{B}  \right ) \\
   + C
\left \|  \frac{(1+(\log N(Y_k))^{p+q} }
{k^p} {\bf 1}_{A^c}\right \|_1\, ,
\end{multline*}
for some positive constant $C$, and consequently
\begin{multline}\label{second}
\E ( |X_{k, \overline x}-X_{k, \overline y} |^q)
\le C \int_{\{\log N(g) \ge  k\}} \left ( \log N(g) \right)^q \, \mu(dg) 
\\ + C  \P(\log d(A_{k-1}\cdot \overline x, A_{k-1}\cdot  \overline y)\ge -\ell k) 
\int_G \left (\log N(g) \right)^{q}
\mu(dg)
\\+ C  \int_{\{\log N(g) <  k\}} \frac{\left (\log N(g) \right)^{p+q}}{k^p}
\mu(dg)\, .
\end{multline}
Now, for $q \in (0,p)$ there exist two positive constants $K$ and 
$L$ such that 
\begin{equation}\label{third}
\sum_{k\ge 1}k^{p-q-1}\int_{\{\log N(g) \ge  k\}} \left ( \log N(g) \right)^q
\mu(dg)
\le K \int_G(\log N)^p\, d\mu<\infty \, ,
\end{equation}
and 
\begin{equation} \label{four}
\sum_{k\ge 1}k^{p-q-1}\int_{\{\log N(g)<  k\}} \frac{(\log N(g))^{p+q}}{k^p}
\mu(dg)\le L \int_G(\log N)^p\, d\mu <\infty \, .
\end{equation}
In the case where $q\in [1,p)$, since $p-q-1 \leq p-2$, we infer from \eqref{second}, 
\eqref{third}, \eqref{four} and 
\eqref{complete-cocycle}  that \eqref{decadix} holds. In the case where $q \leq 1$, since 
$p-q-1 \geq p-2$, the condition \eqref{complete-cocycle} implies \eqref{decadix2}. This completes
the proof of Proposition \ref{newprop}. \hfill $\square$

\smallskip

\noindent {\bf Proof of Proposition \ref{newprop2}.}
Let us first prove \eqref{square-series}. 
Using \eqref{trivial-est}, we see that 
\begin{equation*}
\left | \tilde X_{k, \overline x}^2-
 \tilde X_{k, \overline y}^2  \right |
\le 2 (\log N(Y_k)+ |\lambda_\mu|)  |
\sigma (Y_k ,A_{k-1}\cdot \overline x)-\sigma (Y_k ,A_{k-1}\cdot \overline y)|
\end{equation*}
Proceeding as in \eqref{first} and \eqref{second}, we obtain that
\begin{multline*}
\E \left (\left | \tilde X_{k, \overline x}^2-
 \tilde X_{k, \overline y}^2  \right | \right )
\le C \int_{\{\log N(g) \ge  k\}}  \log N(g) (\log N(g) + |\lambda_\mu| ) \, \mu(dg)\\
+ C \P(\log d(A_{k-1}\cdot  \overline x, A_{k-1}\cdot \overline y)\ge -\ell k) \int  
\log N(g) (\log N(g) + |\lambda_\mu| ) \, \mu(dg)\\ 
+ C \int_{\{\log N(g)<  k\}} \frac{(\log N(g))^{p}(\log N(g) + |\lambda_\mu| )}{k^{p-1}}
\mu(dg)\, ,
\end{multline*}
for some positive constant $C$.
We conclude as in Proposition \ref{newprop} (using similar arguments 
as in \eqref{third} and \eqref{four}).

\smallskip

Let us prove \eqref{mixed-series}. Let $2k\ge i>j\ge k$. 
We start from the simple decomposition
\begin{multline*}
\tilde X_{i, \overline x}
 \tilde X_{j, \overline x}- \tilde X_{i, \overline y}
 \tilde X_{j, \overline y}  =
\tilde X_{i, \overline x}\left(\sigma(Y_j,A_{j-1}\cdot \overline x)- \sigma(Y_j,A_{j-1}\cdot \overline y)\right)
\\  +\left( \sigma(Y_i,A_{i-1}
\cdot \overline x)- \sigma(Y_i,A_{i-1}\cdot \overline y)\right)\tilde X_{j, \overline y}
:= W_{i,j}+Z_{i,j}\, .
\end{multline*}
Using \eqref{trivial-est}, \eqref{regularity} (with $\kappa=p$) and independence, and proceeding as in \eqref{first} and \eqref{second}, we obtain that 
\begin{multline*}
\E(|W_{i,j}|)\le \left(|\lambda_\mu| + \|\log N\|_{1,\mu} \right) \int_{\{\log N(g) \ge  j\}} \log N(g) \, \mu(dg)
\\ +\|\log N\|_{1,\mu}\left(|\lambda_\mu| + \|\log N\|_{1,\mu} \right) \P\left (\log d(A_{j-1}\cdot 
\overline x, A_{j-1}\cdot  \overline y)\ge 
-\ell j/2 \right) \\
+ C\left(|\lambda_\mu| + \|\log N\|_{1,\mu} \right)  \int_{\{\log N(g)<
 j\}} \frac{(\log N(g))^{p+1}}{j^p}
\mu(dg)\, ,
\end{multline*}
and 
\begin{multline*}
\E_\nu(|Z_{i,j}|)\le \left(|\lambda_\mu| + \|\log N\|_{1,\mu} \right) \int_{\{\log N(g) \ge  k\}} \log N(g) \, \mu(dg)
\\ +\|\log N\|_{1,\mu} \E\left((|\lambda_\mu| + \log(N(Y_j))\, {\bf 1}_{\{\log d(A_{i-1}\cdot 
\overline x, A_{i-1}\cdot  \overline y)\ge -i\ell /2\}}\right ) 
\\ + C  \left(|\lambda_\mu| + \|\log N\|_{1,\mu} \right) \int_{\{
\log N(g)<  k\}} \frac{(\log N(g))^{p+1}}{k^p}
\mu(dg)\, ,
\end{multline*}
for some positive constant $C$.
Let $\gamma<p-3+1/p$. It suffices to prove that 
$$
\sum_{k\ge 1} k^\gamma \max_{k \le j< i\le 2k}\E\left (\log(N(Y_j))\, {\bf 1}_{\{\log d(A_{i-1}\cdot 
\overline x, A_{i-1}\cdot  \overline y)\ge -i\ell /2\}}\right ) 
<\infty\, .
$$
Using the H\"older inequality, it is enough to prove that 
$$
\sum_{k\ge 1} k^\gamma \max_{k\le i\le 2k}\left(\P(\log d(A_{i-1}\cdot \overline x, 
A_{i-1}\cdot  \overline y)\ge -i\ell /2)\right)^{(p-1)/p}
<\infty\, .
$$
Using the H\"older inequality again it suffices to find  $\delta >1$ 
such that 
$$
\sum_{k\ge 1} k^{\delta/(p-1)}k^{\gamma p/(p-1)}
\max_{k\le i\le 2k}\left(\P(\log d(A_{i-1}\cdot \overline x, A_{i-1}\cdot \overline y)\ge -i\ell /2)\right)\, .
$$
By Lemma \ref{lem-comp}, it suffices to find $\delta>1$, such that 
$\delta/(p-1)+\gamma p/(p-1)\le p-2$. in particular, it suffices that 
$(p-2)(p-1)-\gamma p>1$, which holds by assumption. 
\hfill $\square$

\section{General results under projective conditions} \label{proj}

In this section, we state general results under projective conditions, that will be needed to prove versions of Theorems \ref{main-theo} and \ref{rates} in stationary regime. Proposition \ref{MZMW} 
is new and is somewhat optimal. Proposition \ref{ASIP} slightly improves previous results. Proposition \ref{DMR} is taken from 
Dedecker, Merlev\`ede and Rio \cite{DMR}. 
Finally Proposition \ref{vbe} is  a new moment  inequality, in the spirit of von Bahr and Esseen \cite{BE},
that will be useful to prove that the results hold for any starting points. 

\smallskip

The proofs of Propositions \ref{MZMW},  \ref{ASIP}
and \ref{vbe}  are given in Section \ref{Sec:inter}.

\smallskip

We shall state  Propositions \ref{MZMW} and \ref{ASIP} in presence of an invertible measure 
preserving transformation, since Proposition \ref{DMR} has been 
proved in that situation. This will be enough for our purpose. 

\smallskip

 Let $(\Omega,\F,\P)$ be a probability space and $\theta$ be an invertible measure preserving transformation. Let $\G_0\subset \F$ be a $\sigma$-algebra, such that $\G_0\subset \theta^{-1}(\G_0)$. For every $n\in \Z$ define 
$\G_n:=\theta^{-n}(\G_0)$.

 For every $Z\in L^p(\Omega,\F,\P)$, we consider the
following maximal functions
\begin{gather}
\label{maxfunp} \M_{p}(Z,\theta):= \sup_{n\ge 1}\frac{\left |\sum_{k=0}^{n-1}
Z\circ \theta^k \right |}{n^{1/p}}\, ,
\quad \mbox{ if $1\le p<2$} .
\end{gather}

Write also $T_n:= Z+\cdots +Z\circ \theta^{n-1}$, and, for any real-valued random variable $V$ and $p\geq 1$, let
$
\|V\|_{p, \infty}= \sup_{t>0} t \left ( {\mathbb P}(|V|>t)\right )^{1/p} 
$.

\smallskip

\begin{prop}\label{MZMW}
Let $1<p < 2$. Let $Z\in L^p(\Omega,\G_0,\P)$ be such that 
\begin{equation}\label{MWp-cond}\sum_{n\ge 1}\frac{\|\E(T_n|\G_0)\|_p}{n^{1+1/p}}<\infty\, .
\end{equation}
There exists a constant $C_{p}>0$, depending only on $p$ such that 
\begin{equation}\label{inemaxLp}
\|\M_p(Z)\|_{p,\infty}\le C_{p} \left (\|Z\|_{p}+
\sum_{n\ge 1}\frac{\|\E(T_n|\G_0)\|_p}{n^{1+1/p}} \right ) \, .
\end{equation}
Moreover, 
\begin{equation}\label{asMZ}
T_n=o\left (n^{1/p} \right) \qquad \mbox{$\P$-a.s.}
\end{equation}
and, there exists $K>0$ such that for every positive integer $d$,
\begin{gather}\label{est11}
\left\|\max_{1\le i\le 2^d} | T_i|\, \right\|_p \le \frac{K}{p-1} 2^{d/p}
\left( \|Z\|_{p}+ \sum_{k=0}^d 2^{-k/p}\|\E(T_{2^k}|\G_{-2^k})\|_{p}\right) \, .
\end{gather}

\end{prop}
\noindent {\bf Remarks.} 
An inequality similar to \eqref{est11} is given in Theorem 3 of \cite{WZ}. 
It is not hard to prove that 
\eqref{MWp-cond} holds as soon as 
\begin{equation}\label{Heydep-cond}
\sum_{n\ge 1}\frac{\|\E(Z\circ \theta^n|\G_0)\|_p}{n^{1/p}}<\infty\, .
\end{equation}
Condition \eqref{MWp-cond} may be seen as an $L^p$-analogue of the so-called 
Maxwell-Woodroofe condition \cite{MW}. As in the papers \cite{PU},
\cite{PUW} or \cite{Cuny-preprint} (see Section D.3), it can be shown 
that \eqref{MWp-cond} is somewhat optimal for \eqref{asMZ}.

\begin{prop}\label{ASIP}
 Let $  2\le p\le 4$ and assume that $\theta$ is ergodic if $p=2$. 
 Let $Z\in L^p(\Omega,\G_0,\P)$ be such that 
\begin{equation}\label{gordin-cond}
\sum_{n\ge 1}\|\E(Z\circ \theta^n|\G_0)\|_p<\infty \, ,
 \quad \text{for $p\in [2,4)$,}
\end{equation} 
and 
\begin{equation}\label{gordin-condbis}
\sum_{n\ge 1}
\log (n) \| \E(Z\circ \theta^n|\G_0)\|_p<\infty \, ,
 \quad \text{for $p=4$.}
 \end{equation}
If $p \in (2,4]$, assume also that
$$\sum_{n\ge 1}\frac{\|\E(T_n^2|\G_0)-\E(T_n^2)\|_{p/2}}{n^{1+2/p}}<\infty
\, .$$ 
Then $\E(T_n^2)/n\longrightarrow \sigma^2$ as 
$n \rightarrow \infty$, and 
\begin{itemize}
\item [$(i)$] If $2\le p< 4$, one can redefine $(T_n)_{n\ge 1}$ without changing its distribution on a (richer) probability space  on which there exists iid random variables  
$(W_n)_{n\ge 1}$ with common distribution ${\mathcal N}(0,\sigma^2)$,  such that  
\begin{equation}\label{asp-eq}
\left |T_n- \sum_{i=1}^n
W_i \right |= o(r_n)
\quad \text{$\P$-a.s.}\, ,
\end{equation}
where $r_n=\sqrt{n\log \log n}$ when  $p=2$ and  $r_n= n^{1/p} \sqrt{\log n}$ 
when $2<p<4$.
\item [$(ii)$] If $p= 4$, one can redefine $(T_n)_{n\ge 1}$ without changing its distribution on a (richer) probability space  on which there exists iid random variables  
$(W_n)_{n\ge 1}$ with common distribution   ${\mathcal N}(0,\sigma^2)$, such that 
\begin{equation} \label{asp-eq2}
\left |T_n- \sum_{i=1}^n 
W_i
\right |= O\left (n^{1/4}\sqrt{\log n} \left(\log \log n \right)^{1/4}  \right)
\quad \text{$\P$-a.s.}\, .
\end{equation}
\end{itemize} 
\end{prop}
\noindent {\bf Remark.}  The condition $
\sum_{n\ge 1}\|\E(Z\circ \theta^n|\G_0)\|_p<\infty$ 
ensures a martingale-coboundary decomposition. It is possible to  weaken 
this condition as done for instance in 
\cite{DDM}. Since in our application the condition $
\sum_{n\ge 1}\|\E(Z\circ \theta^n|\G_0)\|_p<\infty$
 is satisfied, we do not state those refinements.

\smallskip

\begin{prop}\label{DMR}
Let $  2< p\le 3$. Let $Z\in L^p(\Omega,\G_0,\P)$ be such that 
$$\sum_{n\ge 1}\|\E(Z\circ \theta^n|\G_0)\|_p<\infty
\, , \quad \text{for $p \in (2,3)$,}
$$  
  and 
$$
\sum_{n \geq 1} \log (n) \Vert\E (Z\circ \theta^n |{\mathcal G_0})  \Vert_3 < \infty \, , \quad \text{for $p=3$.}
$$
Assume also that 
$$
\sum_{n \geq 1} \frac{\Vert  {\mathbb E}
(T_n^2 | {\mathcal G}_0) - 
{\mathbb E}
(T_n^2 ) \Vert_{p/2}}{n^{3-p/2}}< \infty \, .
$$
Then $n^{-1}\E(T_n^2)\longrightarrow \sigma^2$
as $n \rightarrow \infty$, and, denoting by 
$L_n$ the distribution of $n^{-1/2} T_n$ and by 
$G_{\sigma}$ the normal distribution with mean 
zero and variance $\sigma^2$, one has:
\begin{itemize}
\item  [$(i)$] If $p\in (2,3)$,  then, for any $r \in [p-2, p]$, 
$$
W_r \left( L_n ,G_{\sigma} \right)= O\left (n^{-(p-2)/ 2 \max (1, r)} \right ) \, .
$$
\item [$(ii)$] If $p=3$,  then, for any $r \in (1, 3]$, 
$$
 W_r \left(L_n, G_{\sigma} \right)= O\left (n^{-1/ 2 r} \right ), 
$$
and for $r=1$, 
\begin{equation*}
 W_1 \left(L_n,G_{\sigma} \right)
= O\left (n^{-1/ 2} \log n \right ) \, .
\end{equation*}
\end{itemize}
\end{prop}

To prove Theorem \ref{rates} we shall also need the following von Bahr-Esseen type inequality.
This inequality is stated in the non-starionary case: the $Z_i$'s are real-valued random variables
adapted to an increasing filtration $({\mathcal F}_i)_{i \geq 0}$, and
$T_n=Z_1 + \cdots + Z_n$. 

\begin{prop}\label{vbe}
Let $r \in (1,2]$. The following inequality holds:
$$
\|T_n\|^r_r \leq 2^{2-r} \left ( \sum_{i=1}^n \|Z_i\|^r_r + r \sum_{i=1}^{n-1}
{\mathbb E} \left (|Z_i|^{r-1}|{\mathbb E}(T_n-T_i|{\mathcal F}_i)|\right ) \right ) \, .
$$
Moreover, letting 
$T_n^*= \max (0, T_1, \ldots, T_n)$, 
$$
\|T_n^*\|^r_r \leq \frac 4 {r-1}
\sum_{i=1}^n \|Z_i\|^r_r  + \frac {6r} {r-1} \sum_{i=1}^{n-1}
{\mathbb E} \left (|Z_i|^{r-1}|{\mathbb E}(T_n-T_i|{\mathcal F}_i)|\right ) \, .
$$
\end{prop}

\section{On the convergence of series $\sum_n
n^{-(1+\beta)} \left \|\E(T_n^2|\G_0) 
-\E(T_n^2) \right \|_{p/2}$}

We  keep the same notations as in  previous section. For simplicity, 
if $Z$ belongs to  $L^1(\Omega,\F,\P)$, we shall write 
$Z_n:=Z\circ \theta^n$. 

\smallskip

We want to find conditions relying on series of the type considered in Proposition \ref{newprop}  such that the above series converges for 
a given $p>2$ and a given $\beta\in [1/2,1)$. To do so we shall use 
computations as well as notations from Dedecker, Doukhan and Merlev\`ede 
\cite{DDM}.

\smallskip

For every $k,m\in \N$, define 
\begin{align*}
\gamma_p(m,k) & :=\|\E(Z_mZ_{m+k}|\G_0)-\E(Z_mZ_{m+k})\|_{p/2}\qquad \qquad
\\ \tilde \gamma_p(m) & :=\sup_{m\le j<i\le 2m}\|\E(Z_iZ_{j}|\G_0)-\E(Z_iZ_{j})\|_{p/2}
\, .
\end{align*}
Notice that in our definition of $\tilde \gamma_p(m)$ we take the 
supremum $\sup_{m\le j<i\le 2m}$ while in \cite{DDM} they use 
$\sup_{i\ge j\ge m}$. 

\smallskip

Let $\gamma\in (0,1)$, be fixed for the moment. 
Proceeding as in  (4.18) in \cite{DDM}, we see that (notice that $[m^\gamma]+1
\le 2[m^\gamma]$, for $m\ge 1$)
\begin{multline*}
\left \|\E(T_n^2|\G_0) -\E(T_n^2) \right \|_{p/2} \le 
\sum_{k=1}^n \gamma_p(k,0)+4 \sum_{m=1}^n [m^\gamma] \tilde \gamma_p(m) + 2 \sum_{m=1}^n 
\sum_{k=[m^\gamma]+1}^{n} \gamma_p(m,k)\, ,
\end{multline*}
with the usual convention that an empty sum equals 0. 
We derive that  the sum $\sum_n
n^{-(1+\beta)} \left \|\E(T_n^2|\G_0) 
-\E(T_n^2) \right \|_{p/2}$ is finite provided that 
the following  conditions hold (recall that $\gamma-\beta>-1$):

\begin{align}
\label{cond0} & \sum_{m\ge 1} m^{-\beta}\gamma_p(m,0)<\infty \, ,\\
\label{cond1} 
&\sum_{n\ge 1}n^{\gamma-\beta}\tilde \gamma_p(n)<\infty \, ,\\
\label{cond-inter} 
&\sum_{n\ge 1}\frac{1}{n^{1+\beta} } 
\sum_{m=1}^n 
\sum_{k=[m^\gamma]+1}^{n} \gamma_p(m,k)
<\infty \, .
\end{align}

For every $m,k\in \N$, using the notation $Z_m^{(0)}:=Z_m-\E(Z_m|\G_0) $, define 
$$
\gamma_p^*(m,k):= \left \|\E\left (Z_m^{(0)}Z_{m+k}^{(0)}\Big |\G_0 \right )-\E\left (Z_m^{(0)}Z_{m+k}^{(0)}\right ) \right\|_{p/2}\, .
$$
Writing $P_1(\cdot):=\E(\cdot |\G_1)-\E_{\nu}(\cdot |\G_0)$, and combining (4.20), (4.23) and (4.24) of \cite{DDM}, we infer that, 
\begin{align*}
\sum_{m=1}^n 
\sum_{k=[m^\gamma]+1}^{n} \gamma_p(m,k) & \le 
\sum_{m=1}^n \sum_{k=[m^\gamma]+1}^{n} \sum_{\ell=1}^m \|P_1Z_\ell\|_p\, 
\|P_1Z_{\ell+k}\|_p
\\& \quad \quad + \sum_{m=1}^n \sum_{k=0}^{n}
\|\E(Z_m|\G_0)\|_p\, \|\E(Z_{m+k}|\G_0)\|_p\\
  & \le \sum_{m=1}^n \sum_{k=[m^\gamma]+1}^{n} \sum_{\ell=1}^m \|P_1Z_\ell\|_p\, 
\|P_1Z_{\ell+k}\|_p +\left( \sum_{m=1}^{2n} \|\E(Z_m|\G_0)\|_p\right)^2\, .
\end{align*}
Hence \eqref{cond-inter} holds provided that the following conditions are satisfied
\begin{align}
\label{cond2} &\sum_{n\ge 1}\frac1{n^{1+\beta}}\left( \sum_{m=1}^{2n} 
\|\E(Z_m|\G_0)\|_p\right)^2<\infty\, ,\\
\label{cond-inter2} &\sum_{n\ge 1}\frac1{n^{1+\beta}}\sum_{m=1}^n \sum_{k=[m^\gamma]+1}^{n} \sum_{\ell=1}^m \|P_1Z_\ell\|_p\, 
\|P_1Z_{\ell+k}\|_p<\infty \, .
\end{align}
Now, using that $\sum_{k=[m^\gamma]+1}^{n}\le 
\sum_{k\ge [m^\gamma]+1}$ in the second equation, we see that
\begin{multline*}
\sum_{n\ge 1}\sum_{m=1}^n \sum_{k=[m^\gamma]+1}^{n} \sum_{\ell=1}^m 
n^{-1-\beta}\|P_1Z_\ell\|_p\, \|P_1Z_{\ell+k}\|_p  \\  \le 
C\sum_{m\ge 1} \sum_{k=[m^\gamma]+1}^m \sum_{\ell=1}^m 
m^{-\beta}\|P_1Z_\ell\|_p\, \|P_1Z_{\ell+k}\|_p  \\
 +C \sum_{m\ge 1} \sum_{k\ge m+1} \sum_{\ell=1}^m 
k^{-\beta}\|P_1Z_\ell\|_p\, \|P_1Z_{\ell+k}\|_p \, ,
\end{multline*}
so that 
\begin{multline*}
\sum_{n\ge 1}\sum_{m=1}^n \sum_{k=[m^\gamma]+1}^{n} \sum_{\ell=1}^m 
n^{-1-\beta}\|P_1Z_\ell\|_p\, \|P_1Z_{\ell+k}\|_p
\\
\le C \sum_{k\ge 1}\sum_{\ell\ge 1}\sum_{m=1}^{[k^{1/\gamma}]} 
m^{-\beta}\|P_1Z_\ell\|_p\, \|P_1Z_{\ell+k}\|_p +C \sum_{k\ge 1}\sum_{\ell\ge 1}\sum_{m=1}^{k} 
k^{-\beta}\|P_1Z_\ell\|_p\, \|P_1Z_{\ell+k}\|_p\\
\le \tilde C \sum_{k,\ell\ge 1} k^{(1-\beta)/\gamma}
\|P_1Z_\ell\|_p\, \|P_1Z_{\ell+k}\|_p\le \sum_{\ell\ge 1} 
\|P_1Z_\ell\|_p \sum_{k\ge 1} k^{(1-\beta)/\gamma}\|P_1Z_{k}\|_p\, .
\end{multline*}
Hence, \eqref{cond-inter2} holds as soon as 
\begin{multline*}
\sum_{k\ge 1} k^{(1-\beta)/\gamma}\|P_1Z_{k}\|_p\le C 
\sum_{\ell\ge 0} 2^{\ell((1-\beta)/\gamma}\sum_{ k=2^{\ell}}
^{2^{\ell+1}-1}\|P_1Z_{\ell}
\|_p\\ \le C
\sum_{\ell\ge 0} 2^{\ell((1-\beta)/\gamma+1-1/p} 
\left( \sum_{ k=2^{\ell}}
^{2^{\ell+1}-1}\|P_1Z_{\ell}
\|_p^p\right)^{1/p}<\infty\, ,
\end{multline*}
where we used H\"older for the last inequality.
Applying Lemma 5.2 of \cite{DDM} with $q=p$ we infer that \eqref{cond-inter2} holds as soon as 
 \begin{equation}\label{cond3}
\sum_{n\ge 1} n^{(1-\beta)/\gamma-1/p}\left(\sum_{k\ge n} \frac{\|
\E(Z_k|\G_0)\|_p^p}{k}\right)^{1/p}\, .
\end{equation}
By stationarity, the sequence $(\|
\E(Z_k|\G_0)\|_p)_{k\ge 1}$ is non increasing. Hence, 
using that $\|\cdot\|_{\ell^p}\le \|\cdot \|_{\ell ^1}$, we see that 
\eqref{cond3} holds provided that 
\begin{multline*}
\sum_{\ell \ge 0} 2^{\ell((1-\beta)/\gamma-1/p+1)}\left(\sum_{k\ge \ell} \|
\E(Z_{2^k}|\G_0)\|_p^p\right)^{1/p}  \\
\le \sum_{\ell \ge 0} 2^{\ell((1-\beta)/\gamma-1/p+1)}\sum_{k\ge \ell} 
\|\E(Z_{2^k}|\G_0)\|_p<\infty\, .
\end{multline*}
Changing the order of summation and using again that $(\|
\E(Z_k|\G_0)\|_p)_{k\ge 1}$ is non increasing, we infer that 
\eqref{cond3} holds provided that  
\begin{equation}\label{cond4}
\sum_{k\ge 1} k^{(1-\beta)/\gamma-1/p} \|\E(Z_{2^k}|\G_0)\|_p<\infty\, . 
\end{equation}
Collecting all the above estimates and taking care of Proposition 
\ref{newprop}, we obtain the following result:

\begin{prop}\label{prop} Let $p>2$ and $\beta\in [1/2,1)$. Assume that \eqref{cond0} and \eqref{cond2} 
hold and that there exists $\gamma\in (0,1)$ such that 
\eqref{cond1} and \eqref{cond4}  hold. Then, 
\begin{equation}\label{main-cond}
\sum_{n>0} \frac{\left \|\E(T_n^2|\G_0) 
-\E(T_n^2) \right \|_{p/2}}{n^{1+\beta}}<\infty \, .
\end{equation}
\end{prop}








\section{Proofs of Theorems \ref{main-theo}
and \ref{rates}}

\subsection{Proof of Theorem \ref{main-theo}}\label{Subsec:main-theo}
We first prove a version in stationary regime, i.e. under $\P_\nu$.
The proof makes use of Proposition \ref{MZMW} and Proposition \ref{ASIP}. 
Those results are stated in the context of an invertible dynamical system. 
Let us explain how to circumvent that technical matter. Theorem 
\ref{main-theo} is a limit theorem for the process $(X_n)_{n\ge 1}$, which is a functional of the Markov chain $((Y_n,W_{n-1}))_{n\ge 1}$ with state space $G\times X$ and stationary distribution $\mu\otimes \nu$. 
Since that Markov chain is stationary, it is well-known that, by 
Kolmogorov's theorem, there exists a probability $\hat \P$ on the measurable space $(\hat \Omega,\hat \F)=((G\times X)^{\Z},(\B(G)\otimes \B(X))^{\otimes \Z})$, invariant by the shift $\hat \eta$ on $\hat \Omega$, and 
such that the law of the coordinate process $(\hat V_n)_{n\in \Z}$ (with values in $G\times X$ ) under $\hat \P$ is the
same as the one of the process $((Y_n,W_{n-1}))_{n\ge 1}$ 
under $\P_\nu$. In particular they both are Markov chains. Moreover, 
$(\hat\Omega, \hat \F,\hat\P,\hat \eta)$ is ergodic, which is not difficult to prove.

\smallskip

For every $n\in \Z$, define $\hat X_n:=
 \sigma (\hat V_0 ) \circ 
\hat\eta ^n-\hat \E  (\sigma  (\hat V_0
 ) )$ and $\hat \G_n:=\sigma\{\hat X_k\, :\, k\le n\}$. Then, using the Markov property one can prove easily that for every $p\ge 1$, 
and every $m\ge n\ge 1$,
\begin{align}
\label{change1}\|\hat \E(\hat X_n|\hat \G_0)\|_p&= \| \E_\nu(\tilde X_n |\tilde X_0)
\|_p \le \sup_{\overline x,\overline y\in X} \E
\left (\left |X_{n,\overline x}
-X_{n,\overline y} \right | \right)\, ,\\
\label{change2}
\left \|\hat \E(\hat X_nX_m|\hat \G_0)-\hat \E(\hat X_nX_m) \right \|_p & = \left \| \E_\nu(\tilde X_n\tilde X_m|\tilde X_0)-\E_\nu(\tilde X_n\tilde X_m) )
\right \|_p \\ &  \le \sup_{\overline x,\overline y\in X} \E
\left(\left |\tilde X_{n,\overline x}\tilde X_{m,\overline x}
-\tilde X_{n,\overline y}\tilde X_{m,\overline y}\right |\right)\,.
\nonumber
\end{align}

Let us prove $(i)$. 
Let us apply Proposition \ref{MZMW} with $Z:=\hat X_1$. Notice that by  \eqref{change1}
and Proposition \ref{newprop} (see the remark after it), 
\eqref{Heydep-cond} holds. It follows that 
$$
\hat X_1+\cdots + \hat X_n =o\left (n^{1/p} \right ) \qquad 
\mbox{$\hat \P$-a.s.}
$$
Then, we infer that 
$$
S_n-n\lambda_\mu =o\left (n^{1/p} \right ) \qquad 
\mbox{$\P_\nu$-a.s.}
$$
 or equivalently that for $\nu$-almost every $\overline x
 \in X$,
  $$
S_{n, \overline x}-n\lambda_\mu =o\left ( n^{1/p}
\right ) \quad \mbox{${\mathbb P}$-a.s.}
 $$
 In particular, there exists $y\in \R^d$ with $\|y\|=1$ such that 
  $$
\log \|A_n y\|-n\lambda_\mu =o\left (n^{1/p} \right) \quad \mbox{${\mathbb P}$-a.s.}
 $$

Let $x\in \R^d$ be such that $\|x\|=1$. By Proposition 3.2 page 52 of \cite{BL},  there exists 
a random variable $C$ satisfying  $C(\omega)>0$ for 
${\mathbb P}$-almost every $\omega \in \Omega$,
 and such that, for every $n\in \N$, 
\begin{equation}\label{crucial}
C \le \frac{\| A_n x\|}{\|A_n\|}\le 1 \, .
\end{equation}

Applying this inequality  with $x=y$, we infer that 
$$
\log \|A_n \|-n\lambda_\mu =o\left(n^{1/p}\right) \quad \mbox{${\mathbb P}$-a.s.}
 $$
 and then that for every $x\in \R^d$ such that
 $\|x\|=1$, 
$$
\log \|A_n x \|-n\lambda_\mu =o\left (n^{1/p}\right ) \quad \mbox{${\mathbb P}$-a.s.}
 $$

 Let us prove items $(iii)$ and $(iv)$. Let us apply Proposition \ref{ASIP}
 with $Z=\hat X_1$. Then clearly, the conclusion 
 of Proposition \ref{ASIP} will hold for $Z=\tilde X_1$ and, arguing as above, items $(iii)$ and $(iv)$ of Theorem \ref{main-theo} will follow from Lemma 4.1 
 of Berkes, Liu and Wu \cite{BLW}. 
 
  Notice first that by \eqref{change1} and \eqref{simple-series}, \eqref{gordin-cond} (or  \eqref{gordin-condbis}) holds.  Hence, it remains to check 
 \eqref{main-cond} with $\beta=2/p$, which follows from  Proposition \ref{corollary} below.
 
\begin{prop}\label{corollary}
Let $p>2$ and $\beta\in [1/2,1)$. Take $Z:= \tilde X_1$. 
Then, \eqref{main-cond} holds if $\beta> 3-p$. For instance, 
one may take $\beta=2-p/2$ when $2<p\le 3$ and $\beta=2/p$ when $2<p\le 4$.
\end{prop}
\noindent {\bf Proof of Proposition \ref{corollary}.} Let $\beta>3-p$ (with $\beta\ge 1/2$).  Using \eqref{simple-series} and \eqref{change1}, 
we see that \eqref{cond2} is satisfied, since $\beta>0$. Using, 
\eqref{simple-series}, \eqref{square-series} and \eqref{mixed-series} 
combined with \eqref{change1} and \eqref{change2},
we infer that \eqref{main-cond} holds if 
\begin{gather}
\label{eq1}
-\beta\le p-3\, ;\\ 
\label{eq2}\gamma-\beta\le p-3 +1/p\, ; \\
\label{eq3} (1-\beta)/\gamma -1/p\le p-2\, . 
\end{gather}
Now, \eqref{eq1} holds by assumption and then \eqref{eq2} holds  with
$\gamma=1/p$. It is then not difficult to prove that 
\eqref{eq3} also holds with $\gamma=1/p$. \hfill $\square$

\smallskip

 Let us prove item $(ii)$. By Proposition \ref{newprop}, we have 
 $$
 \sum_{n\ge 1}\left( \int_X \left  |\E_u (X_n)-\lambda_\mu \right |^2\nu(du)\right)^{1/2} < \infty
 \,  .
 $$
 It is well-known then that Gordin's method applies, i.e. that we have a martingale-coboundary decomposition (with respect to $\P_\nu$), and the martingale has stationary and ergodic increments. 
 Hence we have the weak invariance principle under $\P_\nu$ (see \cite{H}) meaning that, for any continuous and bounded function $\varphi$ from $C([0,1])$ to ${\mathbb R}$,  
\begin{equation}\label{wip1}
 \lim_{n \rightarrow \infty} \left | 
 {\mathbb E}_\nu( \varphi(B_{n})) - 
 \int \varphi ( \sigma \varpi) w(d\varpi)
 \right  | = 0 \, ,
\end{equation}
where $w$ is the distribution of a standard Wiener process.

Assume now that $(ii)$ does not hold. Then, there
exists a continuous and bounded function $\varphi_0$
from $C([0,1])$ to ${\mathbb R}$, 
and a sequence $\overline x_n$ of elements of $X$
such that
\begin{equation}\label{wip2}
  \left | 
 {\mathbb E}( \varphi_0(B_{n, {\overline x}_n}))- 
 \int \varphi_0 ( \sigma \varpi) w(d\varpi)
 \right  |  \quad   \text{does not converge to 0 as $n \rightarrow \infty$}.
\end{equation}

Now, if $\psi$ is any bounded and Lipschitz function from $C([0,1])$ to ${\mathbb R}$,
it follows frow the first assertion of Lemma 
\ref{startingpoint} below that 
\begin{equation}\label{wip3}
\lim_{n \rightarrow \infty} 
 \left | 
 {\mathbb E}( \psi (B_{n, {\overline x}_n})) - 
 {\mathbb E}_\nu( \psi(B_n))  \right | =0 \, .
\end{equation}
Putting together \eqref{wip1} and \eqref{wip3}, 
we infer that $B_{n, {\overline x}_n}$ converges in distribution to $\sigma W$, where $W$ is a standard 
Wiener process. This is in contradiction with \eqref{wip2}, 
which completes the proof of $(ii)$.   \hfill $\square$

\smallskip

It remains to prove the  following lemma (note that the first assertion
has already been  proved in \cite{Jan} when $p>2$).

\begin{lem}\label{startingpoint}
Assume that $\mu$ has a moment of order $p \geq 2$.
Then
\begin{equation*}
\sup_{x, y, \|x\|=\|y\|=1}\left \|\log \|A_n x\|-\log \|A_n y\| 
   \right \|_1 < \infty 
 \, ,
\end{equation*}
for $r \in (1,2]$, 
$$
\sup_{x, y, \|x\|=\|y\|=1}\left \|\log \|A_n x\|-\log \|A_n y\| 
   \right \|_r  
   =
   \begin{cases}
   O(1)  \quad  \quad \quad \quad \quad \, \text{if $r \leq p-1$}\\
   O \left( n^{(r+1-p)/r} \right)
   \quad  \text{if $ r > p-1 $},
  \end{cases}
$$
and for $p \in [2,3]$, 
$$
\sup_{x, y, \|x\|=\|y\|=1}\left \|\log \|A_n x\|-\log \|A_n y\| 
   \right \|_p  = 
   O\left( n^{1/p} \right)  \, .
$$
\end{lem}
\noindent{ \bf Proof of Lemma \ref{startingpoint}.}
For any $x, y \in {\mathbb R}^d$  such that 
$\|x\|=\|y\|=1$, 
 one has
$$
\log \|A_n x\|-\log \|A_n y\| = \sum_{k=1}^n 
X_{k,\overline x}-X_{k, \overline y} \, .
$$
Hence 
\begin{equation}\label{b1}
\left \|\log \|A_n x\|-\log \|A_n y\| 
   \right \|_1 \leq 
   \sum_{k=1}^n 
\|X_{k,\overline x}-X_{k, \overline y} \|_1 \, .
\end{equation}
Using \eqref{b1} and \eqref{decadix} (with $q=1$
and $p\geq 2$), the first assertion of Lemma \ref{startingpoint}
follows.

\smallskip

For the case $r \in (1,2]$, we apply Proposition \ref{vbe}. 
Let $s_n(\overline x, \overline y)= \sum_{k=1}^n X_{k,\overline x}-
X_{k, \overline y}$. Then
\begin{multline}\label{normer}
\left \|\log \|A_n x\|-\log \|A_n y\| 
   \right \|_r ^r \leq 2
   \sum_{k=1}^n 
\|X_{k, \overline x}-X_{k, \overline y}\|_r^r\\
+ 4
\sum_{k=1}^{n-1}
\||X_{k, \overline x}-X_{k, \overline y}|^{r-1}\E(s_n(\overline x, \overline y)-s_k( \overline x, \overline y)
|{\mathcal F}_{k})\|_{1}
\, .
\end{multline}
From equality (3.9) in \cite{BQ} (which can also be deduced from \eqref{simple-series}) we infer  that
\begin{equation}\label{mart}
X_{k, \overline x}-X_{k, \overline y}=d_k(\overline x,
\overline y) +\psi(A_{k-1} \overline x, A_{k-1}
\overline y)
- \psi(A_{k}  \overline x, A_{k} \overline y) \, ,
\end{equation}
where $d_k(\overline x, \overline y)$ is ${\mathcal F}_k$-measurable 
and such that ${\mathbb E}(d_k(\overline x, \overline y)|{\mathcal F}_{k-1})=0$, and $\psi$  is a  bounded function 
(with $|\psi | < M$). 
In particular,  it follows  from \eqref{mart} that
$$
\||X_{k,\overline x}-X_{k, \overline y}|^{r-1}\E(s_n(\overline x, \overline y)-s_k(\overline x, \overline y)
|{\mathcal F}_{k})\|_{1} \leq 2M 
\||X_{k, \overline x}-X_{k, \overline y}|^{r-1} \|_{1} \, ,
$$
so that, by \eqref{normer}, 
\begin{equation}\label{norme2bis}
\left \|\log \|A_n  x\|-\log \|A_n   y\| 
   \right \|_r ^r \leq 
   D \sum_{k=1}^n \left ( \|X_{k, \overline x}-X_{k, \overline y} \|_r^r 
   + \||X_{k, \overline x}-X_{k, \overline y}|^{r-1} \|_{1} \right ) \, ,
\end{equation}
for some positive constant $D$.
Applying \eqref{decadix} (with $p\geq 2$ and $q=r$) and   \eqref{decadix2} (with $p\geq2$ and $q=r-1$), we infer that 
$$
\sum_{k=1}^n  
\left ( \|X_{k,\overline x}-X_{k, \overline y} \|_r^r +\||X_{k, \overline x}-X_{k, \overline y}|^{r-1} \|_{1}
\right )
= 
O\left (\max\left (1,n^{(r+1-p)} \right ) \right)\, ,
$$
and the second  assertion of Lemma \ref{startingpoint}
follows from \eqref{norme2bis}.

Let us prove the last assertion.  Let 
$$
Z_{n, x,
 y}= X_{n, \overline x}-X_{n, \overline y} \quad \text{and} \quad 
 T_n( x,  y)=\sum_{k=1}^n Z_{k,  x,  y}:=g_n(x,y, Y_1,  \ldots, Y_n)\, . $$
  With these notations, let
$\hat T_n( x,  y)=g_n(x,y, Y_2, \ldots, Y_{n+1})$. 
Now, it is easy to see that $T_{n}( x,  y)=Z_{1, x,  y} + \hat T_{n-1}(Y_1  x, Y_1  y)$.  Letting
$\psi_p(t)= |t|^p$, we have
$$
 |T_{n}(x,y)|^p=  \left |\hat T_{n-1}(Y_1 x, Y_1 y)\right |^p +Z_{1,x,y} \int_0^1 \psi'_p\left (\hat T_{n-1}(Y_1 x, Y_1 y)+t Z_{1,x,y} \right ) dt \, .
$$
Hence 
\begin{multline*}
|T_{n}(x,y)|^p  \leq 
\left  |\hat T_{n-1}(Y_1x, Y_1 y) \right |^p   + 2^{p-2}  |Z_{1,x,y}|^p 
 +
p 2^{p-2} |Z_{1,x,y}|  \left |\hat T_{n-1}(Y_1 x, Y_1 y)\right |^{p-1} \, .
\end{multline*}
Let $G_{n, p}(x,y)={\mathbb E}\left (|T_n(x,y)|^p\right)$.
Taking the conditional expectation with respect to 
$Y_1$, we get
\begin{equation*}
 \E \left ( |T_{n}(x,y)|^p | Y_1\right ) \leq 
G_{n-1,p}(Y_1x, Y_1 y)   + 2^{p-2}   |Z_{1,x,y}|^p 
 +
p 2^{p-2}  |Z_{1,x,y}|  G_{n-1,p-1}(Y_1x, Y_1 y) \, .
\end{equation*}
Let $u_n= \sup_{x,y, x\neq 0, y \neq 0} G_{n,p}(x,  y)$
and $v_n= \sup_{x,y, x\neq 0, y \neq 0} G_{n,p-1}(x,  y)$.
It follows that 
\begin{equation*}
 \E \left ( |T_{n}(x,y)|^p | Y_1\right ) \leq 
u_{n-1}  + 2^{p-2}   |Z_{1,x,y}|^p 
 +
p 2^{p-2}  |Z_{1,x,y}|  v_{n-1}\, .
\end{equation*}
Taking first the expectation, and then the maximum, we get
\begin{equation*}
u_n \leq u_{n-1} + 2^{p-2} \sup_{x, y, x \neq 0, y \neq 0} \E \left ( |Z_{1,x,y}|^p \right) 
 +
p 2^{p-2} \sup_{x, y, x \neq 0, y \neq 0} \E \left (|Z_{1,x,y}|\right) v_{n-1}\, .
\end{equation*}
Since $p-1 \in [1,2]$, we know from the second assertion of the lemma that 
$
 v_n = O(1) \, . 
$
Consequently, there exists a positive constant $C$ such that
$$
 u_n \leq u_{n-1} + C \, .
$$
 It follows that $u_n=O(n)$, which is the desired result, since 
 $$
 u_n= \sup_{x, y, x \neq 0, y \neq 0}
 \left \|\sum_{k=1}^n 
X_{k,\overline x}-X_{k, \overline y} \right \|_p^p =
\sup_{x, y, \|x\|=\|y\|=1}
\left \|\log \|A_n x\|-\log \|A_n y\| \right \|_p^p \, .
\qquad \square
 $$

\subsection{Proof of Theorem \ref{rates}}
Let $\nu_n$ be the distribution of $n^{-1/2}
(S_n-n \lambda_\mu)$ under ${\mathbb P}_\nu$. As in the proof of Theorem 
\ref{main-theo}, it is enough to apply Proposition 
\ref{DMR} with $Y=\hat X_1$. From 
Proposition \ref{corollary} (with $\beta=2-p/2$) combined with 
\eqref{change1} and \eqref{change2}, we see that the 
assumptions of Proposition \ref{DMR} are satisfied. It 
follows that:
\begin{itemize}
\item  [$(i)$] If $\mu$ has a moment of order $p\in (2,3)$, then, for any $r \in [p-2, p]$, 
$$
 W_r \left( \nu_{n},G_{\sigma} \right)= O\left (n^{-(p-2)/ 2 \max (1, r)} \right ) \, .
$$
\item [$(ii)$] If $\mu$ has a moment of order 3,  then, for any $r \in (1, 3]$, 
$$
 W_r \left(\nu_{n}, G_{\sigma} \right)= O\left (n^{-1/ 2 r} \right ), 
$$
and for $r=1$, 
\begin{equation*}
 W_1 \left(\nu_{n},G_{\sigma} \right)
= O\left (n^{-1/ 2} \log n \right ) \, .
\end{equation*}
\end{itemize}

Recall that $W_0^*$ is 
an element of ${\mathbb R}^d$ such that 
$\overline{ W_0^*}= W_0$ and $\|W_0^*\|=1$.
To prove the results for any starting point, we use the 
following elementary inequalities:

For $r\leq 1$, 
$$
 \sup_{\overline x \in X}  W_r \left(\nu_{n},  
 \nu_{n, \overline x}\right)
   \leq n^{-r/2} \sup_{x, \|x\|=1} \left \|\log \|A_n x\|-\log \|A_n W_0^*\| 
   \right \|_{1, \nu}^r \, .
$$

For $r >1$, 
$$
 \sup_{\overline x \in X}  W_r \left(\nu_{n},  \nu_{n,\overline x}\right)
   \leq n^{-1/2}  \sup_{x, \|x\|=1} \left \|\log \|A_n x\|-\log \|A_n W_0^*\| 
   \right \|_{r, \nu} \, .
$$

From the first assertion of  Lemma \ref{startingpoint}, we infer that:
for $r \leq 1$, 
$$
  \sup_{\overline x \in X}   W_r \left(\nu_{n},  \nu_{n,\overline x}\right)=O\left (n^{-r/ 2 } \right )  \, .
$$
This proves Theorem \ref{rates} for $r\in [p-2,1]$, since in that case
$$
  \sup_{\overline x\in X}   W_r \left(\nu_{n},  \nu_{n,\overline x}\right) 
  =O\left (n^{-(p-2)/2} \right ) \, .
$$

From the last assertion of Lemma \ref{startingpoint}, we infer that: for $p \in (2,3]$,
$$
  \sup_{\overline x \in X}  W_p \left(\nu_{n},  \nu_{n,\overline x}\right)
   =O\left (n^{-(p-2)/2p} \right ) \, .
$$
This proves the result for $r=p$. 

It remains to consider the case $r \in (1,p)$. We use the elementary inequality
$$
\left (W_r \left(\nu_{n},  \nu_{n,\overline x}\right) \right )^r  \leq  \left( W_1 \left(\nu_{n},  \nu_{n,\overline x}\right) \right )^{(p-r)/(p-1)} \left(W_p \left(\nu_{n},  \nu_{n,\overline x}\right)\right)^{p(r-1)/(p-1)} \, .
$$
It follows form the preceding upper bounds for $W_1 \left(\nu_{n},  \nu_{n,\overline x}\right)$ and $W_p \left(\nu_{n},  \nu_{n,\overline x}\right)$
that 
$$
  \sup_{\overline x \in X} \left ( W_r \left(\nu_{n},  \nu_{n,\overline x}\right)\right )^r
   =O\left (n^{-(p-2)(p-r)/2(p-1)} n^{-(p-2)(r-1)/2(p-1)}\right ) 
   =O\left (\, n^{-(p-2)/2}\right ) ,
$$
which concludes the proof. \hfill $\square$

\section{Proofs of the intermediate results} \label{Sec:inter}

\subsection{Proof of Lemma \ref{lem-comp}}


We first recall the following notation: for any $x \in {\mathbb R}^d-\{0\}$ and $g$ in $G$, $ g \cdot \bar x = \overline{g \cdot x} $. 

\smallskip

Since $\int_G \, \log N(g)\, \mu(dg)<\infty$, we may define a 
\emph{bounded}  function $F_1$, by setting 
$$
F_1(\bar x, \bar y)=\int_G \log(d(g\cdot \bar x, g\cdot \bar y )/(d(\bar x, \bar y))) \mu(dg) \qquad \forall \bar x, \bar y\in  X , \bar x \neq  \bar y
\, .
$$
Then, we define a cocycle as follows. For 
every $g\in G$ and every $\bar x, \bar y\in X$ with $\bar x\neq \bar y$, set 
$\sigma_1(g, (\bar x, \bar y)):= \log(d(g\cdot \bar x, g\cdot \bar y)/(d(\bar x, \bar y)))
-F_1(\bar x, \bar y)$.

Finally, write 
$$
\log(d(A_n\bar x,A_n \bar y)/d(\bar x, \bar y)) 
=M_n +R_n \, ,
$$
with
$$
R_n=R_n(\bar x, \bar y):= \sum_{k=1}^n F_1(A_{k-1}  \bar x,A_{k-1}\bar y)\, .
$$
and 
$$M_n:=\sum_{k=1}^n \sigma_1(Y_k,(A_{k-1} \bar x,A_{k-1}\bar y))\,,
$$ and  notice that $(M_n)_{n\ge 1}$ is a martingale in $L^p$, since 
$\mu$ has a moment of order $p$.

\smallskip

Using that $d(\bar x, \bar y)\le 1$, the proposition will be proved if we can prove 
that there exists $\ell >0$, such that

\begin{equation}\label{complete-cocycle-2}
\sum_{k\ge 1} k^{p-2}\max_{k\le j\le 2k}\sup_{\bar x, \bar y\in X, \bar x\neq \bar y} \P \left (
 R_j (\bar x,\bar y)\ge -2\ell k \right )<\infty\, ,
\end{equation}
and 
\begin{equation}\label{complete-cocycle-3}
\sum_{k\ge 1} k^{p-2}\max_{k\le j\le 2k}
\sup_{\bar x, \bar y\in X, \bar x\neq  \bar y} \P(
 |M_j(\bar x, \bar y) |\ge \ell k)<\infty\, .
\end{equation}


\noindent {\bf Proof of \eqref{complete-cocycle-2}.} 
Let $K>0$ be such that $|F_1|\le K$. Let $n\ge 1$ be an integer. Then $|R_n|\le 2nK$ and 
using that $|{\rm e}^x-1-x|\le x^2{\rm e}^{|x|}$ for every $x\in \R$,  we see that, for every $a>0$, 
$$
|\E ({\rm e}^{aR_n})-1-a\E(R_n)|\le a^2K^2{\rm e}^{aK}  .
$$

By Proposition 6.4 $(ii)$ in \cite{BL}, there exists 
$n_0\in \N$ and $\delta>0$, such that 
$$
\sup_{\bar x, \bar y\in X} \E\left (R_{n_0}(\bar x, \bar y) \right) 
\le -\delta\, .
$$

For this $n_0$, we can find $a_0>0$ small enough such that 
$$
\sup_{\bar x, \bar y\in X}\E\left ({\rm e}^{a_0R_{n_0}(\bar x,
\bar y))}\right )\le 1- a_0 \delta/2:=\rho<1$$

Using that $R_{(k+1)n_0}= R_{kn_0}+R_{n_0}\circ \eta^{kn_0}$ and 
conditioning with respect to $\F_{kn_0}$, we infer that 
$$
\sup_{\bar x, \bar y\in X}\E\left ({\rm e}^{a_0R_{(k+1)n_0}(\bar x,
\bar y)} \right)
\le \sup_{\bar x, \bar y\in X}\E\left ({\rm e}^{a_0R_{kn_0}(\bar x,
\bar y)} \right )
\sup_{\bar x, \bar y\in X}\E \left ({\rm e}^{a_0R_{n_0}(\bar x,
\bar y)} \right )\le \rho ^{k+1}\, .
$$

Hence, there exists $C>0$, such that for every $n\in \N$,
$$
\sup_{\bar x, \bar y \in X} \E\left ({\rm e}^{a_0R_n(\bar x, \bar y)}\right )
\le  
C\rho^{n/n_0}\, .
$$

Let $k\ge 1$ and $k\le j\le 2k$ and let $\alpha:= |\log \rho|/(2a_0 n_0)$. Then
$$
\P(R_j(\bar x, \bar y)\ge -\alpha  k)\le {\rm e}^{a_0\alpha k} 
\E \left ({\rm e}^{a_0 R_j(\bar x, \bar y)}\right ) \le C {\rm e}^{a_0\alpha k} \rho^{k/n_0} \leq C \rho^{k/(2n_0)}\, ,
$$
and \eqref{complete-cocycle-2} holds with $\ell =\alpha /2$. 
\hfill $\square$

\smallskip

\noindent {\bf Proof of \eqref{complete-cocycle-3}.} 
The proof  makes use of 
a result about complete convergence for martingales that we recall below. 
 This result represents a very small 
sample of the general situations treated by Alsmeyer 
\cite{Alsmeyer}, and later generalized by Hao and Liu \cite{HL}. 

\smallskip

  Recall that a sequence of random variables $(D_n)_{n\ge 1}$ is said to be dominated by a (non negative) random variable $X$, if there exists $C>0$ such that for every $x>0$, $\P(|D_n|>x)\le C\P(X>x)$.

\smallskip

The next theorem follows directly  from Theorem 2.2 of \cite{HL}.

\smallskip

\begin{theo}[Alsmeyer \cite{Alsmeyer}, Hao and Liu  \cite{HL}]
\label{HL}
Let $(D_n)_{n\in \N}$ be a sequence of  $(\F_n)_{n\in \N }$-martingale differences 
dominated by a variable $X$. For every
$q>1$, every $\gamma\in (1,2]$ and every $L\in \N$, there exists  $C>0$,  such that for every $n\ge 1$ and every $\varepsilon >0$, 
\begin{multline}
 \label{complete-conv}\P\left (\max_{1\le k\le n} |D_1+\cdots +D_k|\ge \varepsilon n \right) 
 \le n \P\left(X>   \frac{\varepsilon n}{4(L+1)}\right)
 \\  +\frac{C}{(\varepsilon n)^{q\gamma(L+1)/(q+L)}} 
\left  \| \E(|D_1|^\gamma|\F_0)+\cdots + \E(|D_n|^\gamma|\F_{n-1}) \right \|_q ^{q(L+1)/(q+L)}\, .
\end{multline}
\end{theo}

We apply Theorem \ref{HL} with $D_k:= \sigma_1(Y_k,(A_{k-1}\bar x,A_{k-1}\bar y))$, $X:= 2\log N(Y_1)$, $\gamma=\min(p,2)$ and $q=L$ 
(to be chosen later). Notice that $(D_n)_{n\in \N}$ is dominated by $X$, 
see for instance Lemma 5.3 page 62 of \cite{BL}. 

\medskip

Since $\E(X^p)<\infty$, it is easy to check that for every 
$\delta>0$,
$$
\sum_{n\ge 1}n^{p-2} \P(X>\delta n) <\infty\, .
$$
Moreover, 
$$\left \| \E(|D_1|^\gamma|\F_0)+\cdots + \E(|D_n|^\gamma|\F_{n-1}) \right \|_q
\le 2n \|X\|_p\, .
$$
Hence, the series 
$$
\sum_{n\ge 1} n^{p-2} \P
\left (\max_{1\le k\le n} |D_1+\cdots +D_k|\ge \varepsilon n
\right) 
$$
converges for every $\varepsilon >0$, as soon as 
$$
\sum_{n\ge 1}  \frac{n^{p-2} }{n^{(\gamma-1)(q+1)/2}}<\infty\, ,
$$
which holds provided that 
$
q> 2(p-1)/(\gamma-1)-1
$.
In particular, we infer that \eqref{complete-cocycle-3} holds 
by taking $\varepsilon =\ell$. \hfill $\square$

\subsection{Proof of Proposition \ref{MZMW}}

We first give a maximal inequality in the spirit of Proposition 2 of  \cite{MP}.
The present form is  just Proposition 4.1 of \cite{Cuny-preprint}. 

\begin{prop}\label{inemax}
Let $X\in L^1(\Omega,\G_0,\P)$. For every $k\ge 0$, write 
$u_k:= |\E(T_{2^k}|\G_{-2^k})|$ and $d_k:= 
\E(T_{2^k}|\G_{-2^k})+(\E(T_{2^k}|\G_{-2^k}))\circ \theta^{2^k} 
-\E(T_{2^{k+1}}|\G_{-2^{k+1}})$. Then, 
for every integer $d\ge 0$, we have 
(with the convention $\sum_{k=0}^{-1}=0$)
\begin{multline}
\max_{1\le i\le 2^d} |T_i|    \le \max_{1\le i\le 2^d}
\left|\sum_{\ell=0}^{i-1} (Z-\E(Z|\G_{-1}))\circ \theta^{\ell } 
\right|
+\sum_{k=0}^{d-1}\max_{1\le i\le 2^{d-k-1}}\left|\sum_{\ell=0}^{i-1} 
 d_k\circ \theta^{2^{k+1}\ell}\right| \\ 
 \label{in-max7}  +u_d + \sum_{k=0}^{d-1}
\max_{0\le \ell \le 2^{d-1-k}-1}u_k\circ \theta^{2^{k+1}\ell} \, .
\end{multline}
In particular, there exists $C>0$, such that for every $p\ge 1$,  
\begin{multline}
 \M_p(X,\theta)    \le C \Big( \sum_{k\ge 0} \frac{u_k}{2^{k/p}} 
+ \sum_{k\ge 0} \frac{\big( \M_1(u_k^p, \theta^{2^{k+1}})
\big)^{1/p}}{2^{k/p}} \\
\label{inemax2}    \qquad +\M_p(X-\E_{-1}(X), 
\theta) +\sum_{k\ge 0} \frac{\M_p( d_k, \theta^{2^{k+1}}}{2^{k/p}}\Big)\, .
\end{multline}
\end{prop}

By Hopf's dominated ergodic theorem (see Corollary 2.2 page 6 of 
\cite{Krengel}), for every $f\in L^1(\Omega,\P)$ and every $k\in \N$, 
$$
\|\M_1(f,\theta^{2^{k}})\|_{1,\infty}\le \|f\|_1\, .
$$
Then, \eqref{inemaxLp} follows from Proposition \ref{inemax} combined 
with Proposition 2.1 of \cite{Cuny-bernoulli}. 

\smallskip

Let us prove \eqref{asMZ}. Define $MW_p:= \{Z\in L^p(\Omega, 
\G_0,\P)~:~\|Z\|_{MW_p}<\infty\}$. Then, 
$(MW_p,\|\cdot \|_{MW_p})$ is a Banach space. 

\smallskip

For every $Z\in L^1(\Omega ,\G_0,\P)$ define $QZ =\E_0(Z\circ \theta)$.  
Notice that $Q^n(Z)=\E_0(Z\circ \theta^n)$. Then, 
clearly $Q$ is a contraction of $L^p(\Omega,\G_0)$.
\medskip
Now, we see that 
\begin{gather*}
\|Z\|_{MW_p}=\sum_{n\ge 0} \frac{\|\sum_{k=0}^{2^n-1}Q^k Z\|_{p}}
{2^{n/p}}\, , \mbox{ if $1<p<2$}\, .
\end{gather*}
Hence, in any case, $Q$ is a contraction  on $MW_p$. 

\smallskip

Writing $V_n:=I+\cdots +Q^{n-1}$ and using that $\|V_nV_kZ\|_{p}\le C\min(k\|V_n\|_{p},n\|V_kZ\|_{p})$, we 
see that, for every $Z\in MW_p$, 
\begin{equation}\label{mean}
\frac{\|V_{2^n} Z\|_{MW_p}}{2^n } \le C_p\left( \frac{\|V_{2^n}Z\|_{p}}{2^{n/p}} 
+ \sum_{k\ge n+1} \frac{\|V_{2^k} Z\|_{p}}{2^{k/p}
}\right)\underset{n\to+\infty} 
\longrightarrow 0\, .
\end{equation}
Now, for every $n\ge 1$, taking $m$ such that $2^m\le n<2^{m+1}$, we have 
 $\|V_{n} Z\|_{MW_p}\le C \sum_{k=0}^m \|V_{2^k}Z\|_{MW_p}=o(n)$.

\smallskip 

In particular, we see that $Q$ is mean ergodic on $MW_p$ and has no 
non trivial fixed point (see e.g. Theorem 1.3 p. 73 of \cite{Krengel}), i.e.,  
\begin{equation}\label{meanerg}
MW_p=\overline{(I-Q)MW_p}^{MW_p} \, .
\end{equation}

Now, by \eqref{inemaxLp} and the Banach principle (see Proposition 
C.1 of \cite{Cuny-bernoulli}) it is enough to prove 
\eqref{asMZ} for a set of elements of $MW_p$ that is dense, in particular 
on $(I-Q)MW_p$. So let $Z=(I-Q)Y$, with $Y\in MW_p$. Then, 
$Z=Y\circ \theta -QY+ Y-Y\circ \theta$, is a martingale-coboundary decomposition in $L^p(\Omega,\P)$. Hence \eqref{asMZ} holds 
since $Y\circ \theta^n =o(n^{1/p})$ $\P$-a.s. 
(by the Borel Cantelli-Lemma) and by the Marcinkiewicz-Zygmund strong law of large numbers for martingales with stationary differences in $L^p$. 

\smallskip

It remains to prove \eqref{est11}. We shall apply once more 
 \eqref{in-max7}. The $L_p$-norm of the first two terms may be 
 estimated thanks to Proposition \ref{vbe}. To estimate the $L_p$-norm of the last term in \eqref{in-max7}, we just notice that 
 
 $$
 \max_{0\le \ell \le 2^{d-1-k}-1}u_k\circ \theta^{2^{k+1}\ell}
 \le \left(\sum_{0\le \ell \le 2^{d-1-k}-1}u_k^p\circ \theta^{2^{k+1}\ell}
 \right)^{1/p}\, ,
 $$
 and \eqref{est11} follows. \hfill $\square$


\subsection{Proof of Proposition \ref{ASIP}}

Since $\sum_{n\ge 1}\|\E(Z_n|\G_{1})\|_p<\infty$, we define a variable 
$R$ in $L^p(\Omega,\F,\P)$ by setting
$$
R:= \sum_{n\ge 1}\E(Z_n|\G_{1})\, .
$$
Then, we have 
$$Z_1= R\circ \theta-\E(R\circ \theta |\G_1) +R-R\circ \theta:= D+R-R\circ \theta\, .
$$
Since $R\in L^p(\Omega,\F,\P)$ it is a standard consequence of the 
Borel-Cantelli lemma that $n^{-1/p}R\circ \theta^n \longrightarrow 0$ $\P$-a.s. as $n$ tends to infinity. Hence, it suffices to prove 
\eqref{asp-eq} with $M_n:= D+\dots +D\circ \theta^{n-1}$ in place of 
$S_n$. 

\smallskip

Since $D\in L^p$, it follows from Theorem 2.1 of Shao \cite{Shao} (see 
the proofs of Corollaries 2.5, 2.7 and 2.8 in Cuny and Merlev\`ede 
\cite{CM}) that we only have to prove that: 
\begin{equation}\label{rate-Lp}
\sum_{i=1}^n (\E(D^2|\G_1)-\E(D^2))\circ \theta^{i-1}=o\left (n^{2/p} \right ) 
\quad \text{$\P$-a.s.}\,,
\end{equation}
if $2\le p <4$, and that 
\begin{equation}\label{rate-L4}
\sum_{i=1}^n (\E(D^2|\G_1)-\E(D^2))\circ \theta^{i-1}=O\left (\left (n\log \log n \right )^{1/2}
\right) 
\quad \text{$\P$-a.s.}\,,
\end{equation}
if $p=4$.

\smallskip

When $p=2$, \eqref{rate-Lp} follows from the ergodic theorem. Now, by Proposition \ref{MZMW} (the ergodicity of $\theta$ is not required) and using orthogonality of martingale increments, \eqref{rate-Lp} holds provided that 
$$
\sum_{n\ge 1}\frac{
\left \|\E(M_n^2|\G_0)-\E(M_n^2  )
\right\|_{p/2}}{n^{1+2/p}}
<\infty\, .
$$
Similary, using Theorem 5.2 of \cite{Cuny-preprint}, 
\eqref{rate-L4} holds provided that 
$$
\sum_{n\ge 1}\frac{
\left \|\E(M_n^2|\G_0)-\E(M_n^2)
\right \|_{2}}{n^{3/2}}
<\infty\, .
$$

Since by assumption, for $p \in (2,4]$,
$$
\sum_{n\ge 1}\frac{\left \|\E(T_n^2|\G_0)-\E(T_n^2)
\right \|_{p/2}}{n^{1+2/p}}
<\infty\, ,
$$
it suffices to prove that 
$$
\sum_{n\ge 1}\frac{\left \|\E(M_n^2|\G_0)-\E(T_n^2|\G_0)
\right \|_{p/2}}{n^{1+2/p}}
<\infty\, .
$$
In case $p \in (2,4)$,  this can be done  as to prove (5.38) in \cite{DMR} (see the proof of Theorem 3.1 in \cite{DMR}).
In case $p=4$,  this can be done  as to prove (5.43) in \cite{DMR} (see the proof of Theorem 3.2 in \cite{DMR}).
\hfill $\square$



\subsection{Proof of Proposition \ref{vbe}.}
We proceed as in the proof of Proposition 1 of 
\cite{DR}. For $r \in (1,2]$, let $\psi_r$ be the function
from ${\mathbb R}$ to ${\mathbb R}^+$ defined by $\psi_r(x)= |x|^r$. We start from the following elementary 
decomposition (using the convention $T_0=0$):
\begin{align*}
|T_n|^r= \psi_r(T_n)& = \sum_{i=1}^n\psi_r(T_i)-\psi_r(T_{i-1})
= \sum_{i=1}^n Y_i\int_0^1 \psi'_r(T_{i-1} +t Y_i) \, dt \\
&= \sum_{i=1}^n Y_i\int_0^1 \left (\psi'_r(T_{i-1} +t Y_i)-\psi'_r(T_{i-1})\right ) \, dt
+ \sum_{i=1}^n Y_i \psi'_r(T_{i-1}) \, .
\end{align*}
Consequently,
\begin{align*}
|T_n|^r &=\sum_{i=1}^n Y_i\int_0^1 \left (\psi'_r(T_{i-1} +t Y_i)-\psi'_r(T_{i-1})\right ) \, dt
+ \sum_{i=1}^n Y_i \left ( \sum_{j=1}^{i-1} \psi'_r(T_{j})-\psi'_r(T_{j-1}) \right)  \\
&=
\sum_{i=1}^n Y_i\int_0^1 \left (\psi'_r(T_{i-1} +t Y_i)-\psi'_r(T_{i-1})\right ) \, dt
+ \sum_{i=1}^{n-1} \left ( \psi'_r(T_{i})-\psi'_r(T_{i-1}) \right)  (T_n-T_i).
\end{align*}
Now, it is easy to check that $|\psi'_r(x)-\psi'_r(y)| \leq r 2^{2-r}|x-y|^{r-1}$. Using this simple fact
and taking the conditional expectation, we obtain 
$$
{\mathbb E}\left (|T_n|^r  \right ) \leq 
 2^{2-r} \sum_{i=1}^n  {\mathbb E}( |Y_i|^r) \int_0^1 r t^{r-1} dt 
 + r 2^{2-r} \sum_{i=1}^{n-1} 
{\mathbb E} \left ( |Y_i|^{r-1}   |{\mathbb E}(T_n-T_i|{\mathcal F}_i)| \right)\, ,
$$
and the inequality is proved.

\smallskip

Let us prove the second inequality. We first write that
$$
(T_n^*)^r = \sum_{i=1}^n (T_i^*)^r-(T_{i-1}^*)^r \, .
$$
Note that for $a \geq b \geq 0$, 
$(r-1)(a^r-b^r) \leq r a(a^{r-1}-b^{r-1})$. Hence,
\begin{equation}\label{max1}
(T_n^*)^r \leq \frac{r}{r-1} \sum_{i=1}^n 
 T_i^*\left ( (T_i^*)^{r-1}-(T_{i-1}^*)^{r-1} \right )
 =
 \frac{r}{r-1} \sum_{i=1}^n 
 T_i\left ( (T_i^*)^{r-1}-(T_{i-1}^*)^{r-1} \right ) \, ,
\end{equation}
the last equality being true because 
$(T_i^*)^{r-1}-(T_{i-1}^*)^{r-1}$ is non zero iff 
$T_i^*=T_i$. Now 
\begin{align}
\sum_{i=1}^n 
 T_i\left ( (T_i^*)^{r-1}-(T_{i-1}^*)^{r-1} \right )
 &=
 \sum_{i=1}^n 
 T_i (T_i^*)^{r-1}- T_{i-1} (T_{i-1}^*)^{r-1} 
 - \sum_{i=1}^n Z_i (T_{i-1}^*)^{r-1} \nonumber \\
 & =
 T_n (T_n^*)^{r-1} -
 \sum_{i=1}^n Z_i (T_{i-1}^*)^{r-1} \, .
 \label{max2}
\end{align}
Recall Young's inequality (based on the concavity of the logarithm): for any $a,  b \geq 0$ and any 
$p, q>1$ 
such that $1/p+1/q=1$, 
$$
ab\le \frac{a^p}{p}+\frac{b^q}{q} \, .
$$
Hence, for $x, y \geq 0$
$$
x y^{r-1} \leq \frac{2^{r-1}}{r} x^r + \frac{r-1}{2 r} y^r \, . 
$$
We infer that 
\begin{equation}\label{max3}
\frac{r}{r-1} |T_n| (T_n^*)^{r-1} \leq 
\frac{2^{r-1}}{r-1} |T_n|^r + \frac 1 2  (T_n^*)^{r}\, .
\end{equation}
Combining \eqref{max1}, \eqref{max2} and \eqref{max3},
we get that 
$$
(T_n^*)^r \leq \frac{2^{r}}{r-1} |T_n|^r - 
\frac{2 r}{r-1} \sum_{i=1}^n Z_i (T_{i-1}^*)^{r-1} \, .
$$
Proceeding as for the first inequality, we get that
$$
(T_n^*)^r \leq \frac{2^{r}}{r-1} |T_n|^r 
- \frac{2 r}{r-1} \sum_{i=1}^{n-1} 
\left ((T_{i}^*)^{r-1} -(T_{i-1}^*)^{r-1})\right) (T_n-T_i) \, .
$$

Since, for $x,y \geq 0$, 
$|x^{r-1}-y^{r-1}| \leq |x-y|^{r-1}$,  we finally get 
that 
$$
\E \left ( (T_n^*)^r \right) 
\leq 
\frac{2^{r}}{r-1}  \E \left ( |T_n|^r \right)
+ 
\frac{2 r}{r-1} \sum_{i=1}^{n-1} 
{\mathbb E} \left ( |Y_i|^{r-1}   |{\mathbb E}(T_n-T_i|{\mathcal F}_i)| \right) \, .
$$
Combining this inequality with the first inequality of the 
proposition, 
the result follows.  \hfill $\square$

\section{Extension of the results}
\label{Sec:MNMC}

In this section, we shall first  explain why the results of Section 
\ref{Sec:Results}
still hold for the sequence $(\log\|A_n\|)_{n\ge 1}$ (with obvious changes in the statements). Next, we shall briefly explain how to deal with $( \log |\langle A_n x,y\rangle|)_{n\ge 1}$.

 \subsection{Matrix norm}
 The fact that the statements of Theorem \ref{main-theo} hold 
 for  $\log\|A_n\|$ instead of $S_{n, \bar x}$ is clear 
 from the proof Theorem \ref{main-theo} (cf. Subsection 
 \ref{Subsec:main-theo}). The crucial 
 point here is Inequality \eqref{crucial}. 
 
 \smallskip
 
 Let $\mu_n$ be the distribution of $\log\|A_n\|$.  
 The fact that the statements of Theorem \ref{rates} hold
 for $\mu_n$ instead of $\nu_{n, \bar x}$ requires some 
 explanations. 
 
Let $\tilde \mu_n$ be the distribution of 
$$
\int_X \log S_{n, u} \, \nu (du)  \, .
$$
 The first point to notice is that the statements of Theorem 
 \ref{rates} are valid for $\tilde \mu_n$ instead $\nu_{n, \bar x}$. This can be proved exactly as for the proof of Theorem
 \ref{rates}, by using some easy consequences of Lemma 
 \ref{startingpoint}, such as
 $$
 \sup_{x,  \|x\|=1}\left \|\log \|A_n x\|- \int_X S_{n, u} \, \nu (du) 
   \right \|_p =O\left ( n^{1/p}\right)  \, , 
 $$
 if $\mu$ has a moment of order $p \in [2,3]$.
 
 The next step is to replace $\tilde \mu_n$ by $\mu_n$. 
 To do this, we need to introduce 
 \begin{equation}\label{delta}
 \delta (\bar x, \bar y):=  \frac{|\langle x,y\rangle|}{\|x\| \, \|y\|} \, .
 \end{equation}
 It follows from Proposition 4.5 of \cite{BQ} that if $\mu$ has a 
moment of order $p>1$,  then 
\begin{equation}\label{deltaBQ}
\sup_{v \in  X} \int_X \left |\log (\delta (u,v)) \right |^{p-1}\, \nu(du) <\infty\, .
\end{equation}
Now, from \cite{BL} pages 52-53, we know that
there exists a random variable $V(\omega)$ with values in 
$X$, such that, for any $x \in {\mathbb R}^d$ such that $\|x\|=1$,
$$
0\le \log\|A_n\|- \log\|A_n x\| \leq \left |\log \delta (\bar x, V)
\right | \, .
$$
Integrating this inequality, we get
$$
\left \|\log\|A_n\|- \int_X S_{n, u} \, \nu (du) \right \|_\infty \leq 
\sup_{v \in  X} \int_X \left |\log (\delta (u,v)) \right |\, \nu(du) < \infty
$$
the term on right hand being finite because $\mu$ has a moment
of order 2. The result easily follows.
 
 \subsection{Results without proximality} Proceeding as in the proof 
 of Theorem 4.11 of \cite{BQ} (using their Lemma 4.13) we infer that 
 the results for matrix norm hold without proximality. 
 Then, we see that the results of Theorems \ref{main-theo} and 
 \ref{rates} hold also without proximality, since \eqref{crucial} 
 do not require proximality but only strong irreducibility. 
 
 \subsection{Matrix coefficients} We shall now explain how to derive results for matrix coefficients, 
i.e. for any given $x,y\in\R^d$ with $\|x\|=\|y\|=1$, we study 
the behaviour of $( \log |\langle A_n x,y\rangle|)_{n\ge 1}$. 

\smallskip

We were not able to extend Theorem \ref{rates} to the matrix 
coefficients. We only succeeded to extend Theorem \ref{main-theo}, but 
under a stronger moment assumption (and it does not seem possible to get rid of the proximality 
assumption here). Our 
argument is inspired by \cite{Jan}. 
We shall use the distance $\delta$ defined in \eqref{delta} and 
the  upper bound \eqref{deltaBQ}.

\smallskip 

Let $1<p\le 4$. Assume  that $\mu$ has a moment of order $p+1$. Let
 us explain why the results from Theorem \ref{main-theo} 
may be extended to the matrix coefficients. Actually, using similar arguments as below, one may see that a moment of order 2 is enough to derive item $(ii)$ of Theorem \ref{main-theo} for the matrix coefficients.

\smallskip
Let $x,y\in\R^d$ such that $\|x\|=\|y\|=1$. We have 
$$
\log |\langle A_n x,y\rangle|= \log \|A_n x\| +\log \|y\| 
+ \log \frac{|\langle A_n x,y\rangle|}{\|A_n x\| \|y\|}\, .
$$
The behaviour of $( \log \| A_n x\|)_{n\ge 1}$ is described in Theorem 
\ref{main-theo}. 

\smallskip

It is obvious (using Lemma 4 of \cite{BLW} to deal with items 
$(iii)$ and$(iv)$) that the results of Theorem \ref{main-theo} will 
hold for the matrix coefficients if we can prove that
$$
\Big|\, \log \frac{|\langle A_n x,y\rangle|}{\|A_n x\| \|y\|}
\, \Big|=o(n^{1/p}) \qquad \mbox{$\P$-a.s.}
$$
or, equivalently, that 
$$
\Big|\, \log \delta(A_n \bar x, \bar y)
\, \Big|=o(n^{1/p}) \qquad \mbox{$\P$-a.s.}
$$
(recall that $A_n \bar x = \overline {A_n x}$).
Since $\delta\le 1$, we are back to prove that  for 
every $\varepsilon >0$, $\P$-a.s., we have 
$$
\delta(A_n \bar x, \bar y)\ge {\rm e}^{-\varepsilon n^{1/p}} 
\qquad \mbox{for all $n$ large enough}\, .
$$

Now, it is well-known (see e.g. Definition 4.1 page 55 and $(9)$ page 61 of 
\cite{BL}), that, for any $x', y'$ in ${\mathbb R}^d- \{0\}$, 

\begin{equation}\label{id}
\frac{|\langle x',y'\rangle|}{\|x'\|\, \|y'\|}= 
\sqrt{1- \left(d\left (\overline{x'},\overline{y'}\right)\right )^2}\, .
\end{equation}
Hence 
$$
 \delta^2(A_n\bar x, \bar y) =1-d^2(A_n\bar x, \bar y) \, .
$$
Now, 
\begin{multline*}
1-d^2(A_n\bar x, \bar y)\ge 1-\big( d(A_n\bar x, W_n)+d( W_n,\bar y)\big)^2\big)
\\
\ge 1-d^2(A_n\bar x, W_n)-2d(A_n\bar x, W_n)d(W_n,\bar y)-d^2(W_n,\bar y)\, .
\end{multline*}

Since for every $n\in \N$, $W_n$  has law $\nu$, 
the variables $(\log \delta (W_n,\bar y))_{n\in \N}$ are identically distributed 
in $L^p(\Omega,\F,\P)$. In particular, it follows from the Borel-Cantelli 
lemma that 
\begin{equation}\label{as-est}
|\log \delta (W_n,\bar y)| = o(n^{1/p}) \quad \text{$\P_\nu$-a.s.}
\end{equation}
Hence, for every $\varepsilon >0$, $\P_\nu$-a.s., we have 
$$
\delta(W_n,  \bar y)\ge {\rm e}^{-\varepsilon n^{1/p}} 
\quad \text{for all $n$ large enough}\, ,
$$
and using \eqref{id} again, for  every $\varepsilon>0$, $\P_\nu$-a.s., 
we have 
\begin{equation}\label{pp}
d(W_n,\bar y)\le \sqrt{ 1-{\rm e}^{-2\varepsilon n^{1/p}} }
\qquad \mbox{for all $n$ large enough}\, .
\end{equation}

As in the proof of Lemma \ref{lem-comp}, we may write 
$$
\log \left( \frac{d(A_n \bar x, W_n)}{d(\bar x, W_0)}\right):= M_n +R_n\, ,
$$
where $(M_n)_{n\ge 1}$ is a (centered) martingale with increments dominated by a variable in $L^{p+1}$ and $(R_n)_{n\ge 1}$ is such that there exists 
$\ell >0$ such that 

$$
\sum_{n\ge 1}\P(R_n \ge -\ell n)<\infty\, .
$$

It is well-known, since $p+1\ge 2$, that $(M_n)_{n\ge 1}$ satisfies the strong law of large numbers (actually even the law of the iterated logarithm). Hence, we have, $\P_\nu$-a.s. 

$$
\log d(A_n \bar x, W_n)\le -\ell n/2 \quad 
\text{for every $n$ large enough}\, .
$$

Finally, we infer that, for every $\varepsilon$, $\P_\nu$-a.s., we have 
$$
1-d^2(A_n\bar x, \bar y)\ge 1- {\rm e}^{-\ell n} -{\rm e}^{-\ell n}\sqrt{ 1-{\rm e}^{-2\varepsilon n^{1/p}} } -
( 1-{\rm e}^{-2\varepsilon n^{1/p}} )\ge C 
{\rm e}^{-2\varepsilon n^{1/p}}\, ,
$$
which is exactly what we wanted to prove. \hfill $\square$


\begin{thebibliography}{99}

\bibitem{Alsmeyer} G. Alsmeyer, \emph{Convergence rates in the law of large numbers for martingales}, Stochastic Process. Appl.  36  (1990),  no. 2, 181-194.

\bibitem{BE} B. von Bahr and  C.-G. Esseen, \emph{Inequalities for
 the $r$-th absolute moment of a sum of random variables, $1\le r\le 2$.}
Ann. Math. Statist  36  1965 299-303.

\bibitem{BQ-book} Y. Benoist and J.-F. Quint, \emph{Random walks on reductive groups} manuscript, http://www.math.u-psud.fr/~benoist/prepubli/prepublication.html

\bibitem{BQ} Y. Benoist and J.-F. Quint, \emph{Central limit theorem 
for linear groups}, accepted for publication in Ann. Probab. 

\bibitem{BLW} I. Berkes, W. Liu and W. B. Wu, \emph{Koml\'os-Major-Tusn\'ady approximation under dependence}, Ann. Probab. 42 (2014), no. 2, 794-817.

\bibitem{BL} P. Bougerol and J. Lacroix, Products of random matrices with applications to Schr\"odinger operators. 
Progress in Probability and Statistics, 8. Birkh\"auser Boston, Inc., Boston, MA, 1985

\bibitem{Cuny-bernoulli} C. Cuny, \emph{ASIP for martingales in 2-smooth Banach spaces. 
Applications to stationary processes},  Bernoulli 21 (2015), no 1, 374-400

\bibitem{Cuny-preprint} C. Cuny, \emph{Limit theorems under the 
Maxwell-Woodroofe condition in Banach spaces}, arXiv:1403.0772, accepted for publication in Ann. Probab.

\bibitem{CM}  C. Cuny and F.  Merlev\`ede, \emph{Strong invariance principles with rate for ``reverse'' martingale differences and applications}, J. Theoret. Probab. 28 (2015), no. 1, 137-183.


\bibitem{DDM} J. Dedecker, P. Doukhan and F. Merlev\`ede,
\emph{Rates of convergence in the strong invariance principle under
projective criteria},  Electron. J. Probab. 17 (2012), no. 16, 31 pp.

\bibitem{DMR} J. Dedecker, F. Merlev\`ede and E. Rio, 
\emph{Rates of convergence for minimal distances in the 
central limit theorem under projective criteria},
Electron. J. Probab. 14 (2009), no. 35, 978-1011.

\bibitem{DR} J. Dedecker and E. Rio, 
\emph{On the functional central limit theorem for stationary processes}, Ann. Inst. H. Poincar\'e Probab. Statist. 36 (2000) 1-34.

\bibitem{FK}  H. Furstenberg and H. Kesten, \emph{Products of random matrices}, Ann. Math. Statist.  31  (1960) 457-469. 

\bibitem{GR}  Y. Guivarc'h and A.  Raugi, \emph{Fronti\`ere de Furstenberg, propri\'et\'es de contraction et th\'eor\`emes de convergence}, 
Z. Wahrsch. Verw. Gebiete  69  (1985),  no. 2, 187-242. 

\bibitem{H} C. C. Heyde, \emph{On the central limit theorem 
and iterated logarithm law for stationary processes}, 
Bull. Austral. Math. Soc. 12 (1975), 1-8. 

\bibitem{HL} S. Hao and Q.  Liu, \emph{Convergence rates in the law of large numbers for arrays of martingale differences}, J. Math. Anal. Appl.  417  (2014),  no. 2, 733-773.

\bibitem{Jan} C. Jan, \emph{Vitesse de convergence dans le TCL pour des 
processus associ\'es \`a des syst\`emes dynamiques ou des 
produits de matrices al\'eatoires}, Th\`ese de l'universit\'e de Rennes 1 (2001), thesis number 01REN10073. 

\bibitem{Jan-CRAS} C. Jan,  Vitesse de convergence dans le TCL pour des cha\^{\i}nes de Markov et certains processus associ\'es à des syst\`emes dynamiques, C. R. Acad. Sci. Paris S\'er. I Math.  331  (2000),  no. 5, 395-398.

\bibitem{KMT} J. Koml\'os, P. Major and G. Tusn\'ady,
\emph{An approximation of partial sums of independent RV's, and the sample DF. I.} Z. Warsch. Verw. Gebiete 32
(1975), 111-131.

\bibitem{Krengel} U. Krengel, Ergodic theorems, de Gruyter Studies in Mathematics,
6. Walter de Gruyter \& Co., Berlin, (1985).

\bibitem{Lepage}  E. Le Page, \emph{Th\'eor\`emes limites pour les produits de matrices al\'eatoires}, Probability measures on groups (Oberwolfach, 1981),  pp. 258–303, Lecture Notes in Math., 928, Springer, Berlin-New York, 1982.

\bibitem{MW} M. Maxwell and M. Woodroofe,  \emph{Central limit theorems for additive functionals of Markov chains}, Ann. Probab. 28 (2000), no. 2, 713-724.

\bibitem{MP} F. Merlev\`ede and M. Peligrad, 
\emph{Rosenthal-type inequalities for the maximum of partial sums of stationary processes and examples}, Ann. Probab. 41 (2013), no. 2,   914-960.

\bibitem{PU}  M. Peligrad and S. Utev, \emph{A new maximal inequality and invariance principle for stationary sequences}, Ann. Probab. 33 (2005), no. 2, 798-815.

 \bibitem{PUW}  M. Peligrad, S. Utev, and W. B. Wu, \emph{A maximal 
 Lp-inequality for stationary sequences and its applications}, Proc. Amer. Math. Soc. 135 (2007), no. 2, 541-550.

\bibitem{Rio} E. Rio, \emph{Inequalities for sums of dependent random variables under projective conditions}, J. Theoret. Probab. 22 (2009), no. 1, 146-163.

\bibitem{Shao} Shao, Q.M. \emph{Almost sure invariance principles for mixing sequences of random variables}. Stochastic Process. Appl. 48 (1993),  319-334.

\bibitem{WZ}  W. B. Wu, and Z. Zhao, \emph{Moderate deviations for stationary processes}, Statist. Sinica 18 (2008), no. 2, 769-782.

\end{thebibliography}
\end{document}